\documentclass[11pt,leqno]{amsart}
\usepackage{amsmath}
\usepackage{amssymb}
\usepackage{amscd}
\usepackage{latexsym}
\usepackage{mathtools}
\usepackage{mathrsfs}
\usepackage{mathdots}
\usepackage{mleftright}
\usepackage{tabularx}
\usepackage{blkarray}
\usepackage{a4wide}
\usepackage[usenames]{color}
\usepackage{epsfig}
\usepackage{enumerate}
\usepackage{enumitem}
\usepackage{subfigure}
\usepackage{url}
\usepackage{tikz,tikz-cd}
\usetikzlibrary{decorations.markings}
\usetikzlibrary{calc}
\usepackage[matrix,arrow,tips,curve]{xy}
\usepackage{pb-diagram,pb-xy}
\usepackage{verbatim}
\usepackage{scalerel}
\usepackage[new]{old-arrows}
\usepackage[colorinlistoftodos,backgroundcolor=orange!20,textsize=scriptsize]{todonotes}
\usepackage{microtype}

\definecolor{cadmiumgreen}{rgb}{0, 0.42, 0.24}
\definecolor{darkred}{rgb}{0.85, 0, 0}
\definecolor{aqua}{rgb}{0, 1, 1}
\definecolor{byzant}{rgb}{0.74, 0.2, 0.64}

\usepackage[
colorlinks, 
linkcolor=blue,
urlcolor=blue,
citecolor=cadmiumgreen,
pdfauthor={Omid Amini, Stéphane Gaubert, and Lucas Gierczak}, 
pdftitle={Equality of tropical rank and dimension for semimodules of tropical rational functions, and computational aspects},
pdfstartview={FitV},
bookmarksdepth=3
]{hyperref}
\usepackage{cleveref}

\newtheorem{theorem}{Theorem}[section]

\newtheorem{lemma}[theorem]{Lemma}
\newtheorem{proposition}[theorem]{Proposition}
\newtheorem{question}[theorem]{Question}
\newtheorem{corollary}[theorem]{Corollary}

\theoremstyle{definition}

\newtheorem{remm}[theorem]{Remark}
\newenvironment{remark}
{\pushQED{\qed}\remm}
{\popQED\endremm}
\newtheorem{exx}[theorem]{Example}
\newenvironment{example}
{\pushQED{\qed}\exx}
{\popQED\endexx}

\newtheoremstyle{italicclaim}% name
	{}% space above
	{}% space below
	{\itshape}% body font
	{}% indent
	{\itshape}% head font
	{.}% punctuation after head
	{.5em}% space after head
	{}% head spec

\theoremstyle{italicclaim}

\newtheorem*{propertyD}{Property D}
\newtheorem*{propertyDbar}{Property $\overline{\mathrm{D}}$}
\newcommand{\supp}{\operatorname{supp}}
\newcommand{\TT}{\mathbb{T}}
\newcommand{\mg}{\Gamma}
\newcommand{\Rat}{\operatorname{Rat}}
\newcommand{\Sl}{\mathrm{sl}}
\newcommand{\Amb}{\operatorname{Amb}}
\renewcommand{\div}{\operatorname{div}}
\newcommand{\ord}{\operatorname{ord}}
\newcommand{\R}{\mathbb{R}}
\newcommand{\Z}{\mathbb{Z}}
\newcommand{\troprank}{\ss[-1pt] r_{\mathrm{trop}}}

\newcommand{\st}{\bigm|}
\newcommand{\T}{\mathrm{T}}
\newcommand{\playermin}{\textsf{Min}}
\newcommand{\playermax}{\textsf{Max}}
\newcommand{\eval}{\vartheta}
\newcommand{\cJ}{\mathcal{J}}
\newcommand{\cK}{\mathcal{K}}
\newcommand{\cI}{\mathcal{I}}
\newcommand{\barphi}{\overline{\varphi}}
\newcommand{\barGamma}{\overline{\Gamma}}

\renewcommand{\geq}{\geqslant}
\renewcommand{\leq}{\leqslant}

\numberwithin{equation}{section}

\renewcommand{\setminus}{\smallsetminus}

\RenewDocumentCommand{\ss}{O{0pt} O{0pt} O{.8} m e{_^}}{
	#4%
	\IfValueT{#5}{
		\sb{\hspace{#1}\scaleobj{#3}{#5}}
	}
	\IfValueT{#6}{
		\sp{\hspace{#2}\scaleobj{#3}{#6}}
	}
}

\makeatletter
\@namedef{subjclassname@2020}{\textup{2020} Mathematics Subject Classification}
\makeatother

\newcommand{\msc}[1]{\href{http://www.ams.org/msc/msc2020.html?t=&s=#1}{#1}}

\title[Equality of tropical rank and dimension for tropical semimodules]{Equality of tropical rank and dimension for semimodules of tropical rational functions, and computational aspects}

\author{Omid Amini}
\address{CNRS--Laboratoire de mathématiques d'Orsay, Université Paris-Saclay}
\email{\href{omid.amini@universite-paris-saclay.fr}{omid.amini@universite-paris-saclay.fr}}
\author{Stéphane Gaubert}
\address{CNRS--CMAP, \'Ecole polytechnique, Institut polytechnique de Paris}
\email{\href{gaubert@cmap.polytechnique.fr}{gaubert@cmap.polytechnique.fr}}
\author{Lucas Gierczak}
\address{I2M, Université d'Aix-Marseille}
\email{\href{lucas.gierczak-galle@univ-amu.fr}{lucas.gierczak-galle@univ-amu.fr}}

\date{First arXiv version April 2025, Revision March 2026}


\begin{document}

\begin{abstract}
	The tropical rank of a semimodule of rational functions on a metric graph mirrors the concept of rank in linear algebra. Defined in terms of the maximal number of tropically independent elements within the semimodule, this quantity has remained elusive due to the challenges of computing it in practice. We establish that the tropical rank is, in fact, precisely equal to the topological dimension of the semimodule, one more than the dimension of the associated linear system of divisors. This implies that the equality of divisorial and tropical ranks in the definition of tropical linear series is equivalent to the pure dimensionality of the corresponding linear system. We then address the question of computing the tropical rank. In particular, we show that checking whether a given family of tropical rational functions is tropically independent is equivalent to solving a turn-based stochastic mean-payoff game, whereas calculating the tropical rank of a finitely generated semimodule of tropical rational functions is $\mathrm{NP}$-hard. We conclude with several complementary results and questions regarding combinatorial and topological properties of the tropical rank.
\end{abstract}

\maketitle

\setcounter{tocdepth}{1}

\section{Introduction} \label{sec:introduction}
	
	Starting from the pioneering work by Baker--Norine~\cite{baker2007riemann} and the subsequent works on algebraic geometry of tropical curves, tropical methods have been quite successful in the study of the geometry of curves and their moduli spaces. We refer to the survey papers~\cite{BJ16, JP21} for a sample of results. The main results of this note are motivated by these developments.
	
	Let $\mg$ be a metric graph (see Section~\ref{sec:preliminaries} for the precise definition). Denote by $\Rat(\mg)$ the union of the set of piecewise linear functions on $\mg$ with integral slopes and the constant function on $\mg$ with value $\infty$ everywhere. Endowed with the two operations of pointwise minimum, denoted by $\oplus$, and pointwise addition of constants, denoted by $\odot$, $\Rat(\mg)$ becomes a semimodule over the semifield of tropical numbers $\TT = (\R \cup \{\infty\}, \oplus, \odot)$.
	
	Let $D$ be a divisor of degree $d$ on $\mg$. The Riemann--Roch space $R(D) \subset \Rat(\mg)$ associated to $D$ is defined by
	\[
		R(D) = \mleft\{f \in \Rat(\mg) \setminus \{\infty\} \st \div(f) + D \geq 0\mright\} \cup \{\infty\}.
	\]
	Here, for $f \neq \infty$, the divisor of $f$, denoted by $\div(f)$, is given by
	\[ \div(f) = \sum_{x \in \mg} \ord_x(f) \, (x), \]
	with the order of vanishing function defined by negative sum of the slopes of $f$ along unit tangent directions to $\mg$ at $x$,
	\[
		\ord_x(f) = -\sum_{\nu \in \T_x \mg} \Sl_\nu f.
	\]
	Endowed with the operations $\oplus$ and $\odot$, $R(D)$ becomes a semimodule over $\TT$.
	
	Let $M$ be a subsemimodule of $R(D)$. The linear system associated to the pair $(D, M)$, denoted by $|(D, M)|$, is defined by
	\[ 
		|(D, M)| = \mleft\{E = \div(f) + D \st f \in M\mright\}.
	\]
	We naturally view $|(D, M)|$ as a subset of the symmetric product
	\[
		\mg^{(d)} = \mg^d/{\mathfrak S_d},
	\]
	where the symmetric group $\mathfrak S_d$ of degree $d$ acts on $\mg^d$ by permuting the coordinates. In the case $M = R(D)$, we get the complete linear system $|D|$.
	
	The symmetric product $\mg^{(d)}$ has a natural polyhedral structure, see~\cite[\S~2.1]{amini2013reduced}. In the case $M = R(D)$, it was proved in~\cite{haase2012linear} that $R(D)$ is finitely generated, and a polyhedral structure for the complete linear system $|D|$ was defined in \emph{loc.~cit.} We will prove in Section~\ref{sec:preliminaries} that for any finitely generated subsemimodule $M \subseteq R(D)$, the linear system $|(D, M)|$ inherits a polyhedral structure as a subspace of $\mg^{(d)}$.
	
	For a finitely generated subsemimodule $M$ in $R(D)$, we define its \emph{dimension} as the topological dimension of the associated linear system, increased by one:
	\[ \dim(M) = \dim(D, M) = \dim |(D, M)| + 1. \]
	More generally, for any subsemimodule $M \subseteq R(D)$, we define its dimension as the supremum (in fact, the maximum, by finite generation of $R(D)$) of $\dim(N)$ over all finitely generated subsemimodules $N$ of $M$:
	\[
		\dim(M) = \sup \, \mleft\{\dim(N) \, \st \, N \text{ is a finitely generated subsemimodule of } M\mright\}.
	\]
	The dimension is an intrinsic numerical invariant of $M$, meaning that it does not depend on the choice of the divisor, see Proposition~\ref{prop:dimension_subsemimodule}. The notation $\dim(M)$ is thus consistent. By analogy with the finitely generated case, we define, for any subsemimodule $M \subseteq R(D)$, $\dim |(D, M)| = \dim(M) - 1$. Equivalently, this is the supremum of $\dim |(D, N)|$ over finitely generated subsemimodules $N \subseteq M$.
	
	We recall from~\cite{JP14} that a family of functions $f_1, \dots, f_r$ in $\Rat(\mg) \setminus \{\infty\}$ is called \emph{tropically dependent} if there exist real numbers $c_1, \dots, c_r$ such that the minimum in
	\[ \min_{1 \leq i \leq r} \, (f_i(x) + c_i) \]
	is achieved at least twice for every $x \in \mg$. Otherwise, the family is called \emph{tropically independent}.
	
	For a subsemimodule $M$ of $\Rat(\mg)$, we define the \emph{tropical rank} $\troprank(M)$ as the maximum integer $r$ for which there exist $r$ tropically independent elements $f_1, \dots, f_r$ in $M$. Also, we define
	\[
		\troprank|(D, M)| = \troprank(M) - 1.
	\]
	
	Our first result is the following theorem (see Section~\ref{sec:proof_main_theorem} for the proof).
	
	\begin{theorem} \label{thm:tropical_rank_equals_topological_rank}
		For each subsemimodule $M \subseteq R(D)$, we have 
		\[ \troprank(M) = \dim(M) \qquad \text{and} \qquad \troprank|(D, M)| = \dim|(D, M)|. \]
		In particular, we have $\troprank|D| = \dim|D|$.
	\end{theorem}
	
	We discuss related results in the literature. Motivated by applications, a combinatorial theory of linear series was developed in recent works~\cite{AG22, JP22}. A combinatorial linear series of rank $r$ in both of these works is a finitely generated subsemimodule $M \subseteq R(D)$ which verifies $r(D, M) = \troprank|(D, M)|$, with some extra constraints. Here, $r(D, M)$ is the divisorial rank of $(D, M)$, introduced originally in~\cite{baker2007riemann} in the case $M = R(D)$. It is defined as the maximum integer $r \geq -1$ such that for any effective divisor $E$ of degree $r$, there exists an element $f \in M$ with $\div(f) + D - E \geq 0$.
	
	In particular, Jensen and Payne studied in~\cite{JP22} a rigid notion of tropical linear series~\cite[Def.~1.5]{JP22}. They proved that for any such linear series, $r(D, M) = \dim|(D, M)|$ holds whenever $M$ is finitely generated. In fact, the proof of Lemma~4.5, \emph{ibid.}, implies more generally that $r(D, M) \leq \dim|(D, M)|$ with no further assumption on $M$, whereas the proof of Corollary~4.7, \emph{ibid.}, implies that $\dim|(D, M)| \leq \troprank|(D, M)|$ when $M$ is finitely generated. \Cref{thm:tropical_rank_equals_topological_rank} shows that the equality $\dim|(D, M)| = \troprank|(D, M)|$ remains valid for arbitrary submodules of $\Rat(D)$, not just for tropical linear series in the sense of~\cite{JP22} -- this theorem applies in particular to situations in which $r(D, M) < \troprank|(D, M)|$.
	
	We also note that Yoshitomi showed that the divisorial rank of a complete tropical linear series is smaller than a certain dimension defined in lattice-theoretical terms~\cite[Thm.~2.7]{yoshitomi}. He also noticed that the divisorial rank is not invariant under isomorphism (Example 6.5, \emph{ibid.}). In contrast, the tropical rank, and thus the dimension, possess such an invariance.
	
	The first equality in Theorem~\ref{thm:tropical_rank_equals_topological_rank} extends to the case of semimodules of rational functions a characterization proved by Develin, Santos and Sturmfels for subsemimodules of $\mathbb{T}^n$~\cite[Cor.~5.4]{develin2007rank}, which was also observed by Butkovi\v{c} in an early exploration of the notion of tropical rank -- see Corollary~6.2.23 and Theorem~6.2.18 in~\cite{Butkovi2010} and the references therein, especially~\cite{Butkovi1985}.
	
	Using Theorem~\ref{thm:tropical_rank_equals_topological_rank}, we prove the following reformulation of the equality $r(D, M) = \troprank|(D, M)|$.
	
	\begin{theorem} \label{thm:pure_dimensionality}
		Let $M \subseteq R(D)$ be a finitely generated subsemimodule. The following statements are equivalent.
		
		\begin{enumerate}
			\item \label{enumitem:equality_divisorial_tropical_rank} The equality $r(D, M) = \troprank|(D, M)|$ holds.
			
			\item \label{enumitem:pure_dimensionality} The linear system $|(D, M)|$ is of pure dimension $r(D, M)$.
		\end{enumerate}
	\end{theorem}
	
	The proof is given in Section~\ref{sec:proof_pure_dimensionality}. It uses Theorem~\ref{thm:comparison_dimension_rank}, which states that all the maximal faces of the polyhedral complex $|(D, M)|$ have dimension at least $r(D, M)$. After the first version of our preprint appeared on arXiv, we learned that Theorem~\ref{thm:comparison_dimension_rank} had been proved by Dupraz in his master thesis~\cite{Dupraz24}. His result is now included in the recent preprint~\cite{Changetal} (see Theorem~1.5 there), which develops further the work initiated in~\cite{JP22}. We include below our proof of Theorem~\ref{thm:comparison_dimension_rank} as it builds on different principles.

	\smallskip
	
	Then, we deal with computability issues related with the tropical rank. We prove the following result.
	
	\begin{theorem} \label{thm:equiv}
		The following problems are equivalent up to polynomial-time Turing reductions:
		\begin{enumerate}
			\item Checking whether a family of tropical rational functions is tropically independent;
			
			\item Solving a turn-based stochastic mean-payoff game.
		\end{enumerate}
	\end{theorem}
	
	We note that turn-based stochastic mean-payoff games are among the fundamental problems in algorithmic game theory with an unsettled complexity: games of this kind have been known to be in the complexity class $\mathrm{NP} \cap \mathrm{coNP}$ since the work of Condon~\cite{condon}, yet it remains unknown whether they belong to P.
	
	One of the reductions stated in~\Cref{thm:equiv}, from tropical linear independence to games, is established in \Cref{thm:reduction} below. This reduction is based on a non-linear fixed-point approach and extends to the metric graph setting an earlier result by Akian, the second author, and Guterman~\cite{AGGut10}, which reduced the tropical linear independence problem for families of vectors of $\mathbb{T}^n$ to \emph{deterministic} mean-payoff games. At the same time, it yields a refinement of the certificates of tropical independence introduced by Farkas, Jensen, and Payne in \cite[Thm.~1.6]{FJP20} by providing a \emph{metric} (\emph{spectral}) measure of independence in terms of a non-linear eigenvalue. This refinement is formulated in Theorem~\ref{thm:certificate_independence_generalized} and may be of independent interest.
	
	The converse reduction is formulated in Theorem~\ref{thm:reduction_converse}. It extends a result of Grigoriev and Podolskii~\cite[Thm.~7]{podolskii}, showing that deterministic mean-payoff games reduce to tropical linear independence for families of vectors of $\mathbb{T}^n$.
	
	Consequently, \Cref{thm:equiv} shows that passing from vectors of $\mathbb{T}^n$ to tropical rational functions on metric graphs corresponds, in algorithmic complexity terms, to passing from deterministic to stochastic games.
	
	As a corollary of the above results, we obtain that checking tropical linear independence of a family of rational functions is in $\mathrm{NP} \cap \mathrm{coNP}$, hence unlikely to be NP-hard (see Corollary~\ref{cor:decision_probleme_independence_NP_coNP}). In contrast, we show that computing the tropical rank of semimodules of rational functions is $\mathrm{NP}$-hard (see Theorem~\ref{thm:NP_hardness}).
	
	We provide in \Cref{thm:expressive} a geometric interpretation of the proof of \Cref{thm:reduction_converse}, showing that any semilinear set stable by tropical linear combinations can be expressed in terms of tropical linear dependence of rational functions.
	
	Finally, in Section~\ref{sec:complementary}, we present complementary results and open questions regarding the combinatorial and topological properties of semimodules related to the tropical rank.
	
	\subsubsection*{Acknowledgments}
		
		The authors thank Matthew Dupraz for helpful discussions.
	
	OA was partially supported by ANR (AdAnAr project, ANR-24-CE40-6184). SG was partially supported by ANR (ZADyG project, ANR-25-CE48-7058). LG was partially supported by ANR (SINTROP project, ANR-22-CE40-0014).
	
\section{Preliminaries} \label{sec:preliminaries}
	
	In this section, we gather some background on metric graphs and their divisor theory. We refer to the survey paper~\cite{BJ16} and~\cite{amini2013reduced, haase2012linear, JP14} for more details. For each positive integer $n$, we denote by $[n]$ the set of positive integers $i$ satisfying $1 \leq i \leq n$.
	
	\subsection{Linear series on metric graphs}
		
		Let $G = (V, E)$ be a finite connected graph with vertex set $V$ and edge set $E$. Let $\ell \colon E \to \ss \R_{> 0}$ be an edge length function, assigning a positive real number $\ell(e)$ to each edge $e$ of the graph.
		
		To the pair $(G, \ell)$, we associate a metric space $\mg$ as follows. For each edge $e \in E$, we place a closed interval $I_e = [0, \ell(e)]$ of length $\ell(e)$ between the two vertices of $e$. The resulting space inherits a natural quotient topology from the topology on the disjoint union of the intervals $I_e$, identifying endpoints according to the adjacency relations in $G$.
		
		Moreover, this topology is metrizable via the path metric, where the distance between any two points in $\mg$ is defined as the length of the shortest path connecting them.
		
		A metric space $\mg$ obtained in this way is called a \emph{metric graph}, and the pair $(G, \ell)$ is called a \emph{model} of $\mg$. Note that a metric graph that is not a singleton has infinitely many different models.
		
		The group of divisors on $\mg$, denoted by $\mathrm{Div}(\mg)$, is the free abelian group generated by the points of $\mg$. Explicitly, it consists of finite linear combinations of points of $\mg$:
		\[
			\mathrm{Div}(\mg) = \mleft\{\sum_{x \in A \subset \mg} n_x (x) \; \st \; n_x \in \mathbb Z \textrm{ and } A \textrm{ a finite set} \mright\}.
		\]
		Here, we write $(x)$ for the generator corresponding to the point $x \in \mg$. For a divisor $D \in \mathrm{Div}(\mg)$ and $x \in \mg$, the coefficient of $(x)$ in $D$ is denoted by $D(x)$. The \emph{support} of $D$, denoted by $\mathrm{Supp}(D)$, is the set of points $x$ with $D(x) \neq 0$. The \emph{degree} of $D$, denoted by $\mathrm{deg}(D)$, is defined as the sum of its coefficients
		\[ \mathrm{deg}(D) = \sum_{x \in \mg} D(x). \]
		
		The complete linear system $|D|$ is given, using the notation of the introduction, by
		\[
			|D| = |(D, R(D))| = \mleft\{\div(f) + D \st f \in R(D)\mright\}.
		\]
		We have a natural embedding 
		\[
			\eta \colon |D|\hookrightarrow \mg^{(d)}
		\]
		which maps each divisor $E = (p_1) + \dots + (p_d)$ in $|D|$ to the corresponding point in the symmetric product $\mg^{(d)}$, given by
		\[
			\eta(E) = (p_1, \dots, p_d).
		\]
		
		\begin{theorem} \label{thm:polyhedral_structure} 
			Let $D$ be a divisor of degree $d$ and $M \subseteq R(D)$ be a finitely generated subsemimodule. Then, $|(D, M)|$ has the structure of a polyhedral space. Moreover, the embedding $|(D, M)| \subseteq \mg^{(d)}$ is piecewise linear.
		\end{theorem}
		
		\begin{proof}
			By~\cite[Thm.~6]{haase2012linear}, $R(D)$ is finitely generated. Let $g_1, \dots, g_l$ be a set of generators for $R(D)$ and $f_1, \dots, f_m$ be a set of generators for $M$. There exist real numbers $\lambda_{ij}$ for $i \in [m]$ and $j \in [l]$ such that
			\[
				f_i = \min_{j \in [l]} \, (g_j + \lambda_{ij}).
			\]
			Now, consider the map
			\[
				\Phi \colon \R^{m} \to \R^{l}
			\]
			given by
			\[
				\Phi(c_1, \dots, c_m) = \mleft(\min_{i \in [m]} (c_i + \lambda_{ij})\mright)_{j = 1}^l.
			\]
			In other words, $\Phi$ is the piecewise linear map from $\R^{m}$ to $\R^{l}$ given by the tropical matrix multiplication from the right by the $m \times l$ matrix $(\lambda_{ij})_{i,j}$.
			
			Define the map
			\[
				\Psi \colon \R^{l} \to |D|
			\]
			by sending each $(x_1, \dots, x_l)$ to the element $D + \div(f)$ in $|D|$ with
			\[
				f = \min_{j \in [l]} \, (g_j + x_j).
			\]
			It follows from the description of the polyhedral structure on $|D|$ given in~\cite{haase2012linear} that $\Psi$ is a piecewise linear map.
			
			Finally, the natural embedding
			\[
				\eta \colon |D| \hookrightarrow \mg^{(d)}
			\]
			is a piecewise linear map, see~\cite{amini2013reduced}.
			
			The subset $|(D, M)|$ of $\mg^{(d)}$ is precisely the image of the composition of the maps we have constructed. More precisely, the subset $|(D, M)| \subseteq \mg^{(d)}$ is the image of the piecewise linear map $\eta \circ \Psi \circ \Phi$. This completes the proof of both statements in the theorem.
		\end{proof}
		
		Here are two observations for future use.
		
		\begin{proposition} \label{prop:dimension_subsemimodule}
			Let $D$ and $D'$ be divisors and $M$ be a closed subsemimodule of $\Rat(\mg)$ included both in $R(D)$ and $R(D')$. Then, $\dim(D, M) = \dim(D', M)$.
		\end{proposition}
		
		\begin{proof}
			The addition by $D' - D$ provides a homeomorphism (in fact, an isomorphism of polyhedral spaces) from $|(D, M)|$ to $|(D', M)|$, from which the proposition follows.
		\end{proof}
		
		\begin{proposition} \label{prop:tropical_rank_maximum}
			For each subsemimodule $M$ in $\Rat(\mg)$, $\troprank(M)$ is the supremum of $\troprank(N)$ over finitely generated subsemimodules $N \subseteq M$.
		\end{proposition}
		
		\begin{proof}
			The result follows directly from the observation that a tropically independent family $f_1, \dots, f_r$ in $M$ remains independent in the finitely generated subsemimodule $\langle f_1, \dots, f_r \rangle$.
		\end{proof}
		
		We also state the following comparison result.
		
		\begin{proposition} \label{prop:comparison}
			Given a point $x \in \mg$ and a unit tangent vector $\nu$ based at $x$, let $n_\nu$ be the number of distinct slopes along $\nu$ taken by functions $f \in M$. We have the inequalities
			\[ r(D, M) + 1 \leq n_\nu \leq \troprank|(D, M)| + 1. \]
		\end{proposition}
		
		\begin{proof}
			The first inequality follows by choosing an effective divisor $E$ made up of $r(D, M)$ distinct points close to $x$ in the direction of $\nu$ and applying the definition of the divisorial rank. The second inequality comes from the fact that functions which coincide at $x$ but have pairwise distinct slopes along $\nu$ are tropically independent (see~\cite[Rem.~6.6]{AG22}).
		\end{proof}
	
	\subsection{Certificates of independence} \label{subsec:certificates}
		
		We will need Theorem~\ref{thm:certificate_independence_generalized} below on certificates of independence.
		
		Let $f_1, \dots, f_n$ denote real-valued functions defined on a set $X$ such that for all pairs $i, j \in [n]$, the difference $f_i - f_j$ is bounded on $X$. For example, we can take $X = \Gamma$ and $f_1, \dots, f_n$ elements of $\Rat(\Gamma)$.
		
		We define the map $T = (T_1, \dots, T_n) \colon \R^n \to \R^n$ by setting, for all $c = (c_1, \dots, c_n) \in \R^n$,
		\begin{align}
			T_i(c) = \sup_{x \in X} \mleft(\min_{j \in [n] \setminus\{i\}} (f_j(x) - f_i(x) + c_j) \mright), \qquad \text{for all } i \in [n]. \label{eq:def_reduce}
		\end{align}
		
		The boundedness of the differences $f_i - f_j$ guarantees
		that $T$ takes finite values.
		
		\begin{theorem} \label{thm:certificate_independence_generalized}
			The following assertions are equivalent.
			
			\begin{enumerate}[label=(\arabic*)]
				\item \label{enumitem:it_1} The functions $f_1, \dots, f_n$ are tropically independent.
				
				\item \label{enumitem:it_2} There exists a vector $c \in \R^n$ and a positive real number $\rho$ such that
				\[
					T_i(c) = \rho + c_i, \qquad \text{for all } i \in [n].
				\]
				
				\item \label{enumitem:it_3} There exists a vector $c \in \R^n$ such that
				\[
					T_i(c) > c_i, \qquad \text{for all } i \in [n].
				\]
				
				\item \label{enumitem:it_4} There exist real numbers $c_1, \dots, c_n$ and $x_1, \dots, x_n \in X$ such that, for all $k \in [n]$, the minimum
				\begin{align} \label{eq:minperm}
					\min_{j \in [n]} (f_j(x_k) + c_j)
				\end{align}
				is achieved only for $j = k$.
				
				\item \label{enumitem:it_5} There exist points $x_1, \dots, x_n \in X$ such that the minimum
				\begin{align} \label{eq:allperm}
					\min_{\sigma \in \mathfrak S_n} \sum_{k \in [n]} f_{\sigma(k)}(x_{k})
				\end{align}
				is achieved by a single permutation.
			\end{enumerate}
		\end{theorem}
		
		This theorem unifies several earlier results. The equivalence of~\ref{enumitem:it_1} and~\ref{enumitem:it_4} is~\cite[Thm.~1.6]{FJP20}; the result is stated there for metric graphs, but the proof applies to the more general setting above. The equivalence of~\ref{enumitem:it_1} and~\ref{enumitem:it_3} was proved in~\cite[Thm~4.12]{AGGut10} when $X$ is finite (but allowing the functions to take infinite values); the entries of the subeigenvector $c$ in~\ref{enumitem:it_3} are precisely the real numbers $c_1, \dots, c_n$ in~\ref{enumitem:it_4}. Still when $X$ is finite, the equivalence between~\ref{enumitem:it_1} and~\ref{enumitem:it_5} follows from Proposition~4.1 and Theorem~5.5 of~\cite{develin2007rank}; it was also established in ~\cite[Thm~2.10]{IRowen} and~\cite[Thm~4.12]{AGGut10}. The eigenvalue $\rho$ in~\ref{enumitem:it_2} provides a quantitative measure of independence, see~\Cref{rk:recall_games} below.
		
		To establish~\Cref{thm:certificate_independence_generalized}, we shall use the Hilbert seminorm (also called ``Hopf's oscillation'') $\|\cdot\|_{\mathrm{H}}$, defined on $\R^n$ by setting
		\begin{align*}
			\|x\|_{\mathrm{H}} = \max_i x_i - \min_j x_j, \qquad x \in \R^n.
		\end{align*}
		
		Let $e = (1, \dots, 1)$ be the point in $\R^n$ with coordinates all equal to one. Observe that $\|x\|_{\mathrm{H}} = 2 \cdot \inf \mleft\{\|x + \lambda e\|_\infty \st \lambda \in \R\mright\}$. Therefore, up to a factor of $2$, the norm $\|\cdot\|_{\mathrm{H}}$ agrees with the norm induced by the supremum norm on the quotient space $\R^n/\R e$.
	 	
		\begin{proof}[Proof of Theorem~\ref{thm:certificate_independence_generalized}]
			We first prove \ref{enumitem:it_1} $\Rightarrow$ \ref{enumitem:it_2}. Let
			\[ 
				M^+ \coloneqq \sup_{\substack{(i, j) \in [n]^2 \\ x \in X}} (f_j(x) - f_i(x)) \qquad \text{and} \qquad M^- \coloneqq \inf_{\substack{(i, j) \in [n]^2 \\ x \in X}} (f_j(x) - f_i(x)).
			\]
			We claim that for each $c \in \R^n$, the vector $d = T(c)$ satisfies $\|d\|_{\mathrm{H}} \leq M^+ - M^-$. Indeed, let $k \in [n]$ be such that $c_{k} = \min_{j \in [n]} c_j$. For each $i \in [n]$, we have $d_i \leq M^+ + c_{k}$ and $d_i \geq M^- + c_{k}$, and so $\|d\|_{\mathrm{H}} \leq M^+ - M^-$.
			
			Since $T$ is order-preserving (for the coordinatewise partial order on $\R^n$) and commutes with the addition of a constant, it is nonexpansive (i.e., $1$-Lipschitz) with respect to $\|\cdot\|_\infty$, see~\cite{crandalltartar}. So it induces a quotient map $\varphi \colon \R^n/\R e \to \R^n/\R e$ that is nonexpansive with respect to $\|\cdot\|_{\mathrm{H}}$. Moreover, $\varphi$ preserves the ball of radius $M^+ - M^-$ with respect to $\|\cdot\|_{\mathrm{H}}$ around $0$. Since this ball is a compact convex subset of $\R^n/\R e$, it follows from Brouwer's fixed-point theorem that $\varphi$ has a fixed point. Hence, there exists a real number $\rho$ and a vector $c \in \R^n$ such that $T_i(c) = \rho + c_i$ for all $i \in [n]$.
			
			We show that $\rho$ is positive. For the sake of contradiction, suppose that $\rho \leq 0$. Then, for all $x \in X$, we have
			\[
				c_i \geq c_i + \rho \geq \min_{j \in [n] \setminus \{i\}} (f_j(x) - f_i(x) + c_j)
			\]
			and thus
			\[ f_i(x) + c_i \geq \min_{j \in [n] \setminus\{i\}} (f_j(x) + c_j), \]
			implying that the minimum $\min_{j \in [n]} (f_j(x) + c_j)$ is achieved at least twice, contradicting the assumption in~\ref{enumitem:it_1}. Therefore, $\rho > 0$.
			
			The implication~\ref{enumitem:it_2} $\Rightarrow$ \ref{enumitem:it_3} is trivial. We next show \ref{enumitem:it_3} $\Rightarrow$ \ref{enumitem:it_4}. Suppose that $T_i(c) > c_i$ holds for all $i \in [n]$. Then, for each $i \in[n]$, we can find a point $x_i \in X$ such that
			\[ \min_{j \in [n] \setminus \{i\}}(f_j(x_i) - f_i(x_i) + c_j) > c_i, \]
			and thus $\min_{j \in [n] \setminus \{i\}} (f_j(x_i) + c_j) > f(x_i) + c_i$, as required.
			
			In order to prove \ref{enumitem:it_4} $\Rightarrow$ \ref{enumitem:it_5},
			observe that the set of minimizing permutations in~\eqref{eq:allperm} is unchanged if we replace $f_j$ by $f_j + c_j$. We can thus assume that $c_j = 0$ for all $j \in [n]$. Then, if~\ref{enumitem:it_4} holds, we have that for all $k \in [n]$, the minimum $\min_{j \in [n]} f_j(x_k)$ is achieved only for $j = k$,
			which entails that the minimum in~\eqref{eq:allperm} is achieved only by the identity permutation, giving~\ref{enumitem:it_5}.
			
			We prove the last implication \ref{enumitem:it_5} $\Rightarrow$ \ref{enumitem:it_1}. Permuting the functions $f_1, \dots, f_n$ if necessary, we may assume that the unique permutation giving the minimum in~\eqref{eq:allperm} is the identity. Suppose, for the sake of contradiction, that $f_1, \dots, f_n$ are tropically
			linearly dependent. Then, translating each $f_i$ by
			a constant, we may assume that for all $k \in [n]$, the minimum $\min_{i \in [n]} f_i(x_k)$ is attained at least twice. Consequently, for each $k \in [n]$, there exists an index $\alpha(k) \in [n] \setminus \{k\}$ achieving the minimum in $\min_{i \in [n]} f_i(x_k)$. Fix an initial index $k_0$ and define $k_s = \alpha(k_{s - 1})$ for $s \geq 1$. This produces an infinite sequence in $[n]$, so it must eventually repeat. Hence, it contains a cycle: there exist distinct elements $k_p, \dots, k_{p + q - 1}$ such that $\alpha(k_{p + q - 1}) = k_p$. This cycle defines a permutation $\sigma$, different from the identity, and we have
			\[ \sum_{k \in [n]} f_k(x_k) = \sum_{k \in [n]} f_{\sigma(k)}(x_k), \] 
			which contradicts~\ref{enumitem:it_4}.
		\end{proof}
	 	
		\begin{remark} \label{rk:recall_games}
			The introduction of $T$ is inspired by the approach of \cite[\S~4.2]{AGGut10} in which methods of non-linear fixed-point theory are applied to study the tropical rank of matrices. In particular, the number $\rho$ in~\ref{enumitem:it_1} may be interpreted as an \emph{additive eigenvalue}, and it is unique. It coincides with the \emph{escape rate} of the operator $T$, that is,
			\begin{align} \label{eq:carac_rho}
				\rho = \lim_{k \to \infty} \mleft(k^{-1} \max_{i \in [n]} T_i(T^k(0))\mright) = \lim_{k \to \infty} \mleft(k^{-1} \min_{i \in [n]} T_i(T^k(0))\mright),
			\end{align}
			where $T^k = T \circ \cdots \circ T$ is the $k$-th iterate of $T$, see~\cite{PFT} and the introductory part of~\cite{akian2023} for background. The existence of the additive eigenpair $(c, \rho)$, needed for the implication~\ref{enumitem:it_1} $\Rightarrow$ \ref{enumitem:it_2}, can be alternatively deduced by applying the general non-linear Perron--Frobenius theorem of~\cite{PFT}. It follows from~\cite[Thm.~20]{akian2023} that, for finite $X$, $\rho$ can be interpreted as a quantitative measure of tropical independence: identifying $X = [m]$, each $f_i$ gives a vector in $\R^m$, and $\rho$ is the minimum over all tropical hyperplanes $H$ in $\R^m$ of the maximum distance, in the Hilbert seminorm, of the vectors $f_i$ to $H$.
		\end{remark}

\section{Proof of Theorem~\ref{thm:tropical_rank_equals_topological_rank}} \label{sec:proof_main_theorem}
	
	Since both $\dim(M)$ and $\troprank(M)$ are given by the maximum of $\dim(N)$ and $\troprank(N)$, respectively, over all finitely generated subsemimodules $N$ of $M$ (see Proposition~\ref{prop:tropical_rank_maximum} for the tropical rank), it will be enough to consider only the case where $M$ is finitely generated.
	
	Thus, in what follows, we assume that $M$ is finitely generated.
	
	\subsection{Proof of the inequality $\dim(M) \geq \troprank(M)$}
	
		We will make use of the certificates of independence.
		
		Let $r = \troprank(M)$, and consider a family of tropically independent elements $f_1, \dots, f_r$ in $M$. By the equivalence of~\ref{enumitem:it_1} and~\ref{enumitem:it_4} in Theorem~\ref{thm:certificate_independence_generalized}, there exist points $x_1, \dots, x_r$ and real numbers $c_1, \dots, c_r$ such that for all $i \in[r]$, the minimum over $j$ in
		\[ \min_{j \in [r]} \, (f_j(x_i) + c_j), \]
		is achieved uniquely at $j = i$.
		
		Let $c = (c_1, \dots, c_r) \in \R^{r}$ and define the subset $B_\varepsilon\subset \R^{r}$ as the set of points $p = (p_1, \dots, p_r)$ satisfying 
		\begin{itemize}
			\item $p_1 = c_1$, and
			
			\item $|p_j - c_j| < \varepsilon$ for each $1 < j \leq r$.
		\end{itemize}
		
		We claim that for sufficiently small $\varepsilon > 0$, the elements $f_p$ in $M$, defined for $p \in B_\varepsilon$ by
		\begin{equation} \label{eq:min} 
			f_p(x) = \min_{1 \leq j \leq r} \, (f_j(x) + p_j) \quad \text{for all } x \in \mg,
		\end{equation}
		are all distinct. To prove this, we choose $\varepsilon > 0$ sufficiently small so that, at $x = x_i$, the minimum in \eqref{eq:min} is achieved uniquely at $j = i$ for all $p \in B_\varepsilon$. This ensures that the function $f_p$ uniquely determines the value of $p$.
		
		Since two elements of $M$ which do not differ by a constant correspond to distinct elements in the linear system $|(D, M)|$, we obtain a piecewise linear embedding of $B_\varepsilon$ into $|(D, M)|$, from which we deduce that $\dim|(D, M)| \geq r - 1$. This implies that $\dim(M) \geq r$, establishing the inequality $\dim(M) \geq \troprank(M)$, as required. \qed
	
	\subsection{Proof of the inequality $\dim(M) \leq \troprank(M)$} \label{subsec:proof_second_inequality}
		
		Let $r = \troprank(M)$. Choose a sequence $x_1, x_2, \dots$ forming a dense subset of $\mg$, and let $f_1, \dots, f_m$ be a set of generators of $M$.
		
		For each integer $K \geq 1$, define the evaluation map
		\begin{align*}
			\eval_K \colon M \longrightarrow \TT^K
		\end{align*}
		that sends each function $f\in M$ to its values at $x_1, \dots, x_K$:
		\[
			\eval_K(f) = (f(x_k))_{1 \leq k \leq K}
		\]
		The semimodule 
		\[
			M_K = \eval_K(M) \subset \TT^K
		\]
		is finitely generated.
		
		For each $K \geq 1$, we define a canonical section $\rho_K$ of the projection $\eval_K$, given by
		\begin{align*}
			\rho_K \colon M_K &\longrightarrow M \\
			g &\longmapsto \inf \, \mleft\{f \in M \st \eval_K(f) \geq g\mright\} = \min \, \mleft\{f \in M \st \eval_K(f) \geq g\mright\}.
		\end{align*}
		Since $M$ is finitely generated, the infimum above is a minimum.
		
		The following proposition provides an explicit expression of the section $\rho_K$.
		
		\begin{proposition} \label{prop:formula_section}
			For each $g \in M_K$, we have
			\begin{align} \label{eq:rho}
				\rho_K(g) = \min_{i \in [m]} \mleft(f_i + c_i^K\mright)
			\end{align}
			where 
			\[
				c_i^K= c_i^K(g) = \max_{k \in [K]} (g(x_k) - f_i(x_k)).
			\]
		\end{proposition}
		
		\begin{proof}
			Let $h$ denote the function on the right-hand side of Equation~\eqref{eq:rho}. By the definition of $c_i^K$, we have
			\[
				f_i(x_k) + c_i^K \geq g(x_k) \qquad \text{for all } k \in [K].
			\]
			From this, we deduce that
			\[
				\eval_K\mleft(f_i + c_i^K\mright) \geq g.
			\]
			This implies the inequality $h \geq \rho_K(g)$.
			
			The other way around, since $f_1, \dots, f_m$ is a generating set for $M$, there exist real numbers $p_1, \dots, p_m$ such that
			\[
				\rho_K(g) = \min_{1 \leq i \leq m}(f_i + p_i).
			\]
			Evaluating at $x_k$, we obtain
			\[
				g(x_k) \leq f_i(x_k) + p_i \qquad \text{for all } k \in [K] \text{ and } i \in [m],
			\]
			which implies
			\[
				c^K_i \leq p_i \qquad \text{ for all } i \in [m].
			\]
			Thus, we conclude that $h \leq \rho_K(g)$, proving the proposition.
		\end{proof}
		
		\begin{proposition}
			For each positive integer $K$, the following equalities hold:
			\begin{align*}
				\eval_K \circ \rho_K = \mathrm{Id}_{M_K} &&\eval_K \circ \rho_K \circ \eval_K = \eval_K,
			\end{align*}
			and
			\begin{align*}
				\rho_K \circ \eval_K \leq \mathrm{Id}_M &&\rho_K \circ \eval_K \circ \rho_K = \rho_K.
			\end{align*}
		\end{proposition}
		
		\begin{proof}
			This is straightforward.
		\end{proof}
		
		\begin{proposition} \label{prop:union_projections}
			We have
			\[
				\bigcup_{K \geq 1} \rho_K(M_K) = M.
			\]
		\end{proposition}
		
		\begin{proof}
			The inclusion $\subseteq$ holds by definition. We prove the reverse inclusion $\supseteq$.
			
			Let $f$ be an element of $M$. We show that $f \in \rho_K(M_K)$ for some $K \geq 1$. There exist real numbers $\lambda_i$ such that 
			\[
				f = \min_{i \in [m]} (f_i + \lambda_i).
			\]
			Replacing each $f_i$ with $f_i + \lambda_i$ and removing the unnecessary terms, we may assume without loss of generality that all $\lambda_i$ are zero and that each $f_i$ contributes to the minimum.
			
			This means that for each $i = 1, \dots, m$, there exists a point $y_i$ of $\mg$ such that
			\[
				f_i(y_i) < f_j(y_i) \quad \text{for all } j \neq i. 
			\]
			
			Thus, for each $i$, there exists an open set $U_i \subset \mg$ such that for all $x \in U_i$, we have
			\[
				f_i(x) < \min_{j \neq i} f_j(x).
			\]
			Since the sequence $(x_k)$ is dense in $\mg$, we can choose an index $k_i \geq 1$ such that $x_{k_i} \in U_i$. Let $K = \max_i k_i$.
			
			Now, consider $g = \eval_K(f) \in M_K$. We will show that $c_i^K(g) = 0$ for each $i = 1, \dots m$, from which it follows that $f = \rho_K(g)$.
			
			For a given element $i \in [m]$, we have, for each $k \in [K]$,
			\[
				g(x_k) = f(x_k) \leq f_i(x_k),
			\]
			where equality holds for $k = k_i$. By Proposition~\ref{prop:formula_section}, this implies that $c_i^K(g) = 0$, as required.
		\end{proof}
		
		\begin{proposition} \label{prop:dimension_M_K}
			We have
			\[ \dim(\rho_K(M_K)) \leq r. \]
		\end{proposition}
		
		\begin{proof}
			Recall that $r$ denotes the tropical rank of $M$. Since $M_K$ is a projection of $M$, we have
			\[ \troprank(M_K) \leq \troprank(M). \]
			Thus, to prove the proposition, it suffices to show that
			\[
				\dim(\rho_K(M_K)) \leq \troprank(M_K).
			\]
			The space $M_K$ consists of the column space of the matrix
			\[
				A_K = \Bigl(f_1(x_k), \dots, f_m(x_k)\Bigr)_{1 \leq k \leq K}.
			\]
			By~\cite[Thm.~4.2]{develin2007rank}, $M_K$ is a polyhedral complex of dimension
			\[
				\troprank(A_K) = \troprank(M_K).
			\]
			To conclude, observe that
			\[
				\rho_K \colon M_K \to \rho_K(M_K) \hookrightarrow M
			\]
			is a piecewise linear isomorphism, as given by the explicit formula for $c_i^K(g)$ (see Proposition~\ref{prop:formula_section}). Consequently, 
			\[
				\dim(\rho_K(M_K)) = \dim(M_K),
			\]
			which completes the proof.
		\end{proof}
		
		Using the results stated above, we are now in the position to prove the inequality
		\[
			\dim(M) \leq \troprank(M),
		\]
		finishing the proof of Theorem~\ref{thm:tropical_rank_equals_topological_rank}.
		
		As previously stated, we have an injective piecewise linear map
		\[
			\varphi \colon M/\R \longhookrightarrow \mg^{(d)}, \qquad [f] \longmapsto [\div(f) + D].
		\]
		
		By Proposition~\ref{prop:union_projections}, we obtain
		\[ \varphi(M/\R) = \bigcup_{K \geq 1} \varphi (\rho_K(M_K)/\R). \]
		Moreover, by Proposition~\ref{prop:dimension_M_K}, we have
		\[
			\dim(\varphi(\rho_K(M_K)/\R)) \leq r - 1 \qquad \text{for every $K \geq 1$}.
		\]
		
		\smallskip
		
		Now, let $E$ be a cell of $\varphi(M/\R)$ of maximal dimension. Proceeding by contradiction, assume that $\dim(E) \geq r$. Consider the intersection
		\[
			F_K = \varphi(\rho_K(M_K)/\R) \cap E.
		\]
		Since each $\varphi(\rho_K(M_K)/\R)$ has dimension at most $r - 1$, it follows that $F_K$ has empty relative interior in $E$.
		
		By the Baire category theorem, the countable union $\bigcup_K F_K$ has empty relative interior in $E$. However, this contradicts the fact that $\bigcup_K F_K = E$. Thus, we conclude that $\dim(M) \leq \troprank(M)$, which completes the proof. \qed

\section{Proof of Theorem~\ref{thm:pure_dimensionality}} \label{sec:proof_pure_dimensionality}
	
	Let $D$ be a divisor of degree $d$ on a metric graph $\Gamma$, and $M \subseteq R(D)$ be a finitely generated subsemimodule. By Theorem~\ref{thm:polyhedral_structure}, $|(D, M)|$ has a polyhedral structure. Denote by $r(D, M)$ the divisorial rank of $(D, M)$. We have the following result.
	
	\begin{theorem}[Dupraz, Corollary~2.95.2 in~\cite{Dupraz24}] \label{thm:comparison_dimension_rank}
		Keeping the same notation as above, all the maximal faces of the polyhedral structure on $|(D, M)|$ have dimension at least $r(D, M)$.
	\end{theorem}
	
	We will give an alternate proof below using the Baire category theorem again. Combining this result and Theorem~\ref{thm:tropical_rank_equals_topological_rank}, we deduce Theorem~\ref{thm:pure_dimensionality}.
	
	\begin{proof}[Proof of Theorem~\ref{thm:pure_dimensionality}]
		The implication \eqref{enumitem:pure_dimensionality} $\Rightarrow$ \eqref{enumitem:equality_divisorial_tropical_rank} follows directly from Theorem~\ref{thm:tropical_rank_equals_topological_rank}.
		
		For the reverse implication \eqref{enumitem:equality_divisorial_tropical_rank} $\Rightarrow$ \eqref{enumitem:pure_dimensionality}, note that the equality $\troprank(D, M) = \dim|(D, M)|$, given by Theorem~\ref{thm:tropical_rank_equals_topological_rank} again, implies that $\dim|(D, M)| = r(D, M)$. By applying Theorem~\ref{thm:comparison_dimension_rank}, we conclude that all the maximal faces of $|(D, M)|$ have dimension exactly $r(D, M)$, which implies that $|(D, M)|$ is of pure dimension $r(D, M)$.
	\end{proof}
	
	\begin{proof}[Proof of Theorem~\ref{thm:comparison_dimension_rank}]
		We assume that $|(D, M)|$ is nonempty and $r(D, M) \geq 1$; otherwise, the statement holds trivially. Fix an element $E \in |(D, M)|$, and express it as $E = D + \div(h)$ for $h \in M$. We claim that every closed neighborhood $\mathcal U$ of $E$ in $|(D, M)|$ has dimension at least $r(D, M)$, which will prove the theorem.
		
		For each positive integer $r \leq r(D, M)$, let $S_r$ denote the subset of $\Gamma^{r}$ consisting of all tuples $\underline x = (x_1, \dots, x_r) \in \Gamma^r$ for which there exists a tuple $\underline y = (y_1, \dots, y_{d - r}) \in \Gamma^{d - r}$ such that the sum $(x_1) + \dots + (x_r) + (y_1) + \dots + (y_{d - r})$ belongs to $\mathcal U$. We claim that $S_r$ is a closed subset of $\Gamma^r$.
		
		To see this, suppose we have a sequence $\underline x_t = (x_{1,t}, \dots, x_{r,t})$ of points of $S_r$ converging to $\underline x \in \Gamma^r$. For each $t$, there is a corresponding tuple $\underline y_t = (y_{1,t}, \dots, y_{d - r,t}) \in \Gamma^{d - r}$ such that the sum $(x_{1,t}) + \dots + (x_{r,t}) + (y_{1,t}) + \dots + (y_{d - r,t})$ lies in $ \mathcal U$. By compactness of $\Gamma^{d - r}$, and passing to a subsequence if necessary, we can ensure that $\underline y_t$ converges to some point $\underline y = (y_1, \dots, y_{d - r}) \in \Gamma^{d - r}$. Since $\mathcal U$ is closed, the limit $(x_1) + \dots + (x_r) + (y_1) + \dots + (y_{d - r})$ of the divisors $(x_{1,t}) + \dots + (x_{r,t}) + (y_{1,t}) + \dots + (y_{d - r,t})$ belongs to $\mathcal U$, implying that $\underline x \in S_r$. Therefore, $S_r$ is closed in $\Gamma^r$.
		
		We proceed by induction on $r$ and demonstrate that for each closed neighborhood $\mathcal U$ of an element $E \in |(D, M)|$ as above with corresponding sets $S_j$, $j \leq r(D, M)$, the set $S_r$ contains a nonempty open subset of $\Gamma^r$, i.e., it has a nonempty interior in $\Gamma^r$. By the definition of $S_{r(D, M)}$, this implies that each $\mathcal U$ is of dimension at least $r(D, M)$, as claimed.
		
		We begin by considering the case $r = 1$. Since $r(D, M) \geq 1$, there exists an element $F = D + \div(f)$ in $|(D, M)|$, distinct from $E$, with $f \in M$. For each $t \in \R$, define $f_t = \min \{t + h, f\}$ and let $E_t = D + \div(f_t)$. We observe that $f_t \in M$ and $E_t \in |(D, M)|$.
		
		For $t$ approaching $-\infty$, we have $E_t = E$, and for $t$ approaching $\infty$, we have $E_t = F$. Thus, we obtain a one-dimensional segment in $\mathcal U$, with one endpoint at $E$, completing the proof in this case.
		
		Assume $r \geq 2$ and that the claim has been proved for $r - 1$. Consider the projection map $\pi \colon S_r \to S_{r - 1}$, given by $(x_1, \dots, x_{r - 1}, x_r) \mapsto (x_1, \dots, x_{r - 1})$. For each point $\underline x = (x_1, \dots, x_{r - 1})$ in $S_{r - 1}$, the fiber $\pi^{-1}(\underline x)$ is nonempty. We claim that there exists a nonempty open subset $\mathcal W$ of $\Gamma^{r - 1}$ contained in $S_{r - 1}$ over which each fiber of $\pi$ contains a one-dimensional segment of $\Gamma$ (possibly different for distinct points).
		
		To prove this, let $\mathcal U'$ be a closed neighborhood of $E$ included in the interior of $\mathcal U$. The subset $S'_{r - 1}$ associated to $\mathcal U'$, defined analogously to the previous construction, is a subset of $S_{r - 1}$, and by the induction hypothesis, it has a nonempty interior, which we will show to be the desired nonempty open subset $\mathcal W$. We need to show that for each $\underline x \in S'_{r - 1}$, the fiber $\pi^{-1}(\underline x)$ contains a one-dimensional segment.
		
		Let $\underline y \in \Gamma^{d - r + 1}$ be a point such that $E' = (x_1) + \dots + (x_{r - 1}) + (y_1) + \dots + (y_{d - r + 1}) \in \mathcal U'$. Since $\mathcal U'$ is contained in the interior of $\mathcal U$, $\mathcal U$ is a closed neighborhood of $E'$. Write $E' = D + \div(f)$ for some $f \in M$ and let $M' = M - f$ be the set consisting of all elements of the form $g - f$ for $g \in M$. Then, $M'$ is a finitely generated subsemimodule of $R(E')$.
		
		Moreover, the pair $((y_1) + \dots + (y_{d - r + 1}), M')$ has divisorial rank at least one. Therefore, repeating the argument from the case $r = 1$, there exists a segment of dimension one in $\mathcal U$ of the form $(x_1) + \dots + (x_{r - 1}) + (y_{1,t}) + \cdots + (y_{d - r + 1,t})$ with one endpoint at $E'$. This implies that $\pi^{-1}(\underline x)$ contains a one-dimensional segment.
		
		Recall that $\mathcal W$ is a nonempty open subset of $\Gamma^{r - 1}$ contained in $S_{r - 1}$ over which each fiber of $\pi$ contains a one-dimensional segment of $\Gamma$. Let $(I_j)_{j = 1}^\infty$ be a countable collection of one-dimensional segments in $\Gamma$ such that each one-dimensional segment of $\Gamma$ contains at least one of the $I_j$'s, for some $j \in \mathbb N$. In particular, each fiber of $\pi$ over $\mathcal W$ contains one of the segments $I_j$.
		
		For each $j$, define $\mathcal A_j$ as the subset of $\mathcal W$ consisting of all the points $\underline x$ with $I_j \subseteq \pi^{-1}(\underline x)$. By continuity of the projection map $\pi$ and the closedness of $S_r$, each $\mathcal A_j$ is a closed subset of $\mathcal W$. Moreover, the union of the sets $\mathcal A_j$ covers the full open set $\mathcal W$. Therefore, by the Baire category theorem, there exists some $j \in \mathbb N$ such that $\mathcal A_j$ has a nonempty interior.
		
		We conclude that the product $\mathcal A_j \times I_j$ is contained in $S_r$. Moreover, it has a nonempty interior in $\Gamma^r$, as required.
	\end{proof}

\section{Computational complexity of tropical linear independence and rank} \label{sec:computational_complexity_tropical_linear_independence}
	
	\subsection{From tropical independence to stochastic mean-payoff games}
		
		Checking whether a finite family of tropical vectors is tropically independent reduces to solving a deterministic mean-payoff game, a well-studied example of a repeated game, see~\cite[Thm.~4.12]{AGGut10}. The converse also holds: solving a deterministic mean-payoff game reduces to the problem of checking tropical linear independence~\cite{podolskii}. The question of the existence of a polynomial-time algorithm to solve deterministic mean-payoff games, first raised in~\cite{gurvich}, remains unsettled. When formulated as a decision problem, deterministic mean-payoff games belong to the complexity class $\mathrm{NP} \cap \mathrm{coNP}$~\cite{zwick}, making them unlikely to be $\mathrm{NP}$-complete. Note that tropical vectors with finite entries can be identified to rational functions over a trivial metric graph -- with connected components reduced to isolated vertices. In contrast, the metric graphs we consider here are non-trivial and connected.
		
		The following result shows that for rational functions over metric graphs, checking tropical linear independence reduces to solving a more expressive class of games, with \emph{stochastic} transitions.
		
		\begin{theorem} \label{thm:reduction}
			Checking whether rational functions $f_1, \dots, f_n$ over a metric graph are tropically independent reduces in polynomial time to solving a stochastic turn-based mean-payoff game.
		\end{theorem}
		
		Some technical details on the encoding of $f_1, \dots, f_n$ and the metric graph are in order. We assume that the metric graph is given by a model with rational edge lengths. The input comprises one distinguished vertex. Furthermore, we assume that each of the functions $f_1, \dots, f_n$ is described by a collection of intervals on which it is affine, as well as by the collection of integral slopes on these intervals, together with a rational value at the distinguished vertex. In this way, each function is uniquely determined.
		
		Theorem~\ref{thm:reduction} also leads to effective methods for checking tropical independence, see the discussion in Remark~\ref{rk:rk_compute}.
		
		We now provide some background on the stochastic game appearing in Theorem~\ref{thm:reduction}. It is a special case of the model introduced by Shapley~\cite{shapley_stochastic}. There are two players, called \playermin\ and \playermax. The game is played on a finite state space $[n]$. In each state $i \in [n]$, a finite set of possible actions of \playermax, $A_i$, is specified. Moreover, for each choice of state $i \in [n]$ and for each choice of action $\alpha \in A_i$ of \playermax, a finite set of possible actions of \playermin, $B_{i,\alpha}$, is specified. For every pair of actions $\alpha \in A_i$ and $\beta \in B_{i,\alpha}$, we associate a real number $r_i^{\alpha \beta}$ and a probability measure on the set $[n]$, $P_i^{\alpha \beta} = (P_{i,1}^{\alpha \beta}, \dots, P_{i,n}^{\alpha \beta}) \in \R_{\geq 0}^n$; so we have $\sum_{j \in [n]} P_{i,j}^{\alpha \beta} = 1$.
		
		We consider the \emph{turn-based} game, also known as the \emph{perfect information game}, in which the two players choose their actions sequentially, being informed of the current state and previous actions of the other player. Every stage of the game is played as follows. At a given stage, the current state being $i \in [n]$, \playermax\ selects an action $\alpha \in A_i$, and \playermin, aware of this choice, selects an action $\beta \in B_{i,\alpha}$. Then, \playermin\ pays $r_i^{\alpha \beta}$ to \playermax, and with probability $P_{i,j}^{\alpha \beta}$ the next state becomes $j \in [n]$.
		
		In the finite horizon game, starting from a given initial state, \playermax\ seeks to maximize the expected payment received from \playermin\ over a given number of consecutive stages, whereas \playermin\ seeks to minimize the same payment.
		
		The dynamic programming operator of this game, known as the \emph{Shapley operator}, is the map $T = (T_1, \dots, T_n) \colon \R^n \to \R^n$ defined by
		\begin{align}
			T_i(c) \coloneqq \max_{\alpha \in A_i} \min_{\beta \in B_{i,\alpha}} \bigg(r_i^{\alpha \beta} + \sum_{j \in [n]} P_{i,j}^{\alpha \beta}c_j\bigg).
			\label{eq:def_Shapley}
		\end{align}
		
		The game started from state $i$ and played over $N$ consecutive stages has a value denoted by $v_i^N$, and the value vector $v^N = (v^N_1, \dots, v^N_n) \in \R^n$ satisfies $v^N = T^N(0)$, with $T^N$ being the $N$-th iterate of $T$, see Lem.~VII.1.3 and Rem.~VII.1.1 in~\cite{sorin_repeated_games}. Moreover, the \emph{mean payoff} of this stochastic game is defined as $\chi_i \coloneqq \lim_{N \to \infty} v^N_i/N$. The existence of this limit can be derived as a special case of a result by Bewley and Kohlberg~\cite{bewley_kohlberg}, which applies to concurrent games. In the case of turn-based games, the limit $\chi_i$ is a rational number with a bit size polynomially bounded in the input size.
		
		\begin{remark}
			The mean-payoff problem can be formulated either \emph{quantitatively}, as computing the mean payoff $\chi_i$ of a given initial state $i$ -- or as deciding whether $\chi_i$ exceeds a given threshold --, or \emph{strategically}, as searching for optimal policies for the players. All these variants are equivalent up to polynomial-time (Turing) reductions~\cite{andersson2009complexity}.
		\end{remark}
		
		In order to prove Theorem~\ref{thm:reduction}, we first show that the operator $T$ defined in Equation \eqref{eq:def_reduce} in Section~\ref{subsec:certificates} can be rewritten in the form~\eqref{eq:def_Shapley}. To see this, we consider a finite family $(X_\alpha)_{\alpha \in A}$ of closed segments on edges of $\Gamma$ on which every map $f_j$ is affine and such that $\cup_{\alpha \in A} X_\alpha = \Gamma$. Thus, after identifying the segment $X_\alpha$ with an interval of the form $[u_\alpha, v_\alpha] \subset \R$, the restriction of $f_j$ to the interval $X_{\alpha}$ is of the form $f_j(x) = m_{\alpha,j} x + \eta_{\alpha,j}$ for all $x \in X_\alpha$, for some integer $m_{\alpha,j}$ and some rational number $\eta_{\alpha,j}$. For each $\alpha \in A$, define the operator $T_{\alpha} = (T_{\alpha,1}, \dots, T_{\alpha,n}) \colon \R^n \to \R^n$ by setting
		\begin{align}
			T_{\alpha,i}(c) \coloneqq \sup_{x \in [u_\alpha, v_\alpha]} \mleft(\min_{j \in [n] \setminus \{i\}} \bigl((m_{\alpha,j} - m_{\alpha,i}) x + \eta_{\alpha,j} - \eta_{\alpha,i} + c_j\bigr)\mright),
			\label{eq:part}
		\end{align}
		so that $T_i(c) = \max_{\alpha \in A}T_{\alpha,i}(c)$.
		
		Evaluating the operators $T_{\alpha,i}$, for $\alpha \in A$ and $i \in [n]$ involves solving an optimization problem of the form
		\begin{align} \label{eq:lppb}
			\sup_{x \in [u, v]} \Phi(x) \quad \text{where} \quad \Phi(x) = \min_{j \in J} (\gamma_j x + d_j) \text{ for } x \in [u, v],
		\end{align}
		with $(\gamma_j)_{j \in J}$ and $(d_j)_{j \in J}$ finite families of real numbers, $J \subseteq [n]$, and $u, v \in \R$.
		
		In fact, this optimization problem has an explicit solution given by the following lemma.
		
		\begin{lemma} \label{lem:dual_program}
			Keeping the same notation as above, the value of the optimization problem~\eqref{eq:lppb} is given by
			\begin{align} \label{eq:linear_program}
				\min \mleft(\min_{j \in J^-} (\gamma_j u + d_j), \;
				\min_{j \in J^+} (\gamma_j v + d_j), \;
				\min_{\substack{j \in J^+, k \in J^- \\ \gamma_j > \gamma_k}} \mleft(\frac{-\gamma_k}{\gamma_j - \gamma_k} d_j + \frac{\gamma_j}{\gamma_j - \gamma_k} d_k\mright)\mright)
			\end{align}
			where $J^+ = \mleft\{j \in J \st \gamma_j \geq 0\mright\}$ and $J^- = \mleft\{j \in J \st \gamma_j \leq 0\mright\}$.
		\end{lemma}
		
		\begin{proof}
			The function $\Phi$ is piecewise linear and concave. Hence, the maximum locus of $\Phi$ forms a segment $I$ possibly reduced to a point. If $u \in I$, the minimum in~\eqref{eq:linear_program} is achieved by the first term, which is the value of $\Phi$ at $u$. Similarly, if $v \in I$, the minimum in~\eqref{eq:linear_program} is achieved by the second term, which is the value of $\Phi$ at $v$. It remains to treat the case where $u, v \notin I$. In this case, $I$ is exactly the locus of points of $[u, v]$ around which both outgoing slopes of $\Phi$ are nonpositive. Then, the minimum is achieved by the third term in~\eqref{eq:linear_program}, which is the (constant) value of $\Phi$ on $I$.
		\end{proof}
		
		\begin{remark} \label{rk:dual}
			Lemma~\ref{lem:dual_program} may be interpreted via linear programming: the terms in the minimum correspond to extreme points of the dual of the linear programming formulation of the optimization problem~\eqref{eq:lppb}.
		\end{remark}
		
		We now define a game such that $T_i(c)$ can be rewritten in the form~\eqref{eq:def_Shapley}. The set of states is $[n]$. The set of actions of \playermax\ is $A_i = A$ independently of the state $i \in [n]$. Now, for each $\alpha \in A_i$, to express $T_{\alpha,i}(c)$ in the form~\eqref{eq:lppb}, we take $u = u_\alpha$, $v = v_\alpha$, $J = [n] \setminus \{i\}$, $\gamma_j = m_{\alpha,j} - m_{\alpha,i}$ and $d_j = \eta_{\alpha,j} - \eta_{\alpha,i} + c_j$, and we define $J^+_{\alpha,i} = \mleft\{j \in J \st \gamma_j \geq 0\mright\}$ and	$J^-_{\alpha,i} = \mleft\{j \in J \st \gamma_j \leq 0\mright\}$, as in Lemma~\ref{lem:dual_program}. The set of actions $B_{i, \alpha}$ is the disjoint union of the sets $J^+_{\alpha,i}$, $J^-_{\alpha,i}$ and the set of ordered pairs $(j, k) \in J^+_{\alpha,i} \times J^-_{\alpha,i}$ with $\gamma_j > \gamma_k$. The instantaneous payments and transition probabilities are defined as follows, for all $i \in [n]$ and $\alpha \in A$.
		
		\begin{itemize}
			\item If $\beta = j$ with $j \in J^+_{\alpha,i}$, we set $r_i^{\alpha,\beta} = (m_{\alpha,j} - m_{\alpha,i}) u + \eta_{\alpha,j} - \eta_{\alpha,i}$ and $P_i^{\alpha,\beta} = e_j$, where $e_j$ denotes the $j$-th vector of the canonical basis of $\R^n$.
			
			\item If $\beta = j$ with $j \in J^-_{\alpha,i}$, we set $r_i^{\alpha,\beta} = (m_{\alpha,j} - m_{\alpha,i}) v + \eta_{\alpha,j} - \eta_{\alpha,i}$ and $P_i^{\alpha,\beta} = e_j$.
			
			\item If $\beta = (j, k)$ with $(j, k) \in J^+_{\alpha,i} \times J^-_{\alpha,i}$ and $\gamma_j > \gamma_k$, we define $\pi_j = -\gamma_k/(\gamma_j - \gamma_k)$ and $\pi_k = \gamma_j/(\gamma_j - \gamma_k)$, so that $\pi_j$ and $\pi_k$ are nonnegative and satisfy $\pi_j + \pi_k = 1$. We set $r_i^{\alpha,\beta} = \pi_j (\eta_j - \eta_i) + \pi_k(\eta_k - \eta_i)$ and $P^{\alpha,\beta}_{i} = \pi_je_j + \pi_k e_k$.
		\end{itemize}
		
		As noted in Remark~\ref{rk:recall_games}, the number $\rho$ such that $\rho + c = T(c)$ for some $c \in \R^n$ is unique. Moreover, by~\eqref{eq:carac_rho}, $\rho$ coincides with the value of the turn-based stochastic mean-payoff game, independently of the choice of the initial state $i \in [n]$. By Theorem~\ref{thm:certificate_independence_generalized}, the family $f_1, \dots, f_n$ is tropically linearly independent if, and only if, $\rho$ is strictly positive.
		
		The number of states and actions, as well as the bit sizes of the instantaneous payments and transition probabilities of this game are polynomially bounded in the input size. It follows that verifying whether $f_1, \dots, f_n$ are tropically linearly independent reduces in polynomial time to checking whether the value of a turn-based mean-payoff game is strictly positive. This finishes the proof of Theorem~\ref{thm:reduction}. \qed
		
		\smallskip
		
		We have the following consequence of Theorem~\ref{thm:reduction}.
		
		\begin{corollary} \label{cor:decision_probleme_independence_NP_coNP}
			The decision problem associated with tropical independence of tropical rational functions belongs to the complexity class $\mathrm{NP} \cap \mathrm{coNP}$.
		\end{corollary}
		
		This is a variation on a result of Condon~\cite{condon}, showing that \emph{simple stochastic games}, a subclass of stochastic turn-based games, belong to $\mathrm{NP} \cap \mathrm{coNP}$. We provide a self-contained proof of this corollary, referring the reader to~\cite{andersson2009complexity} for more information on the complexity of turn-based games.
		
		\begin{proof}[Proof of Corollary~\ref{cor:decision_probleme_independence_NP_coNP}]
			We associated to an instance of this decision problem a Shapley operator $T$ of the form~\eqref{eq:def_Shapley}, with the property that there exists a vector $c \in \R^n$ such that $T(c) = \rho + c$, $\rho$ being the (uniquely defined) value of the game. Owing to the piecewise-linear character of $T$ and to the rational character of the parameters $r_{i}^{\alpha,\beta}$ and $P_{i,j}^{\alpha \beta}$ arising in the definition of $T$, the set of solutions $(c, \rho)$ of $T(c) = \rho + c$ can be written as a finite union of rational polyhedra.
			
			Any non-empty rational polyhedron has a \emph{short element} or \emph{short vector}, that is, a rational vector with a bit size polynomially bounded in the size of the external description of the polyhedron. Note that if this polyhedron does not contain an affine line, any vertex provides such a short vector, see e.g.~\cite[Lemma~6.2.4]{Groetschel1993}, and one may always reduce to this case after modding out by the lineality space of the polyhedron, see the proof of Theorem~6.5.7, \emph{ibid.}.
			
			Any such short element $c$ can serve as a certificate of independence: verifying this certificate consists in checking that the vector $T(c) - c$ is constant and strictly positive, which can be done in polynomial time. This shows that the tropical independence problem is in $\mathrm{NP}$. The same vector can serve as a certificate of dependence: we still check that the vector $T(c) - c$ is constant, but now require that this constant be nonnegative. This shows that the tropical independence problem is in co-$\mathrm{NP}$.
		\end{proof}
		
		\begin{remark} \label{rk:rk_compute}
			To solve the game here, it suffices to solve the additive eigenproblem $T(c) = \rho + c$. This can be done using several algorithms that are efficient in practice, see e.g.~\cite{Chatterjee2014} and the discussion in Section~6 of~\cite{akian2023}.
		\end{remark}
		
		\begin{remark}
			\Cref{thm:certificate_independence_generalized} and \Cref{thm:reduction} actually imply a more precise result than \Cref{cor:decision_probleme_independence_NP_coNP}: tropical linear independence belongs to $\mathrm{UP} \cap \mathrm{coUP}$, where the complexity class $\mathrm{UP} \subset \mathrm{NP}$ refers to the languages recognized in polynomial time by a nondeterministic Turing machine with exactly one accepting computation path for every word in the language, and $\mathrm{coUP}$ denotes the complementary class.
			
			We sketch the proof, leaving details to the reader. Instead of considering the certificate $(c, \rho)$ as in the proof of \Cref{cor:decision_probleme_independence_NP_coNP}, we introduce an artificial discount factor $0 < \alpha < 1$, observe that the map $x \mapsto T(\alpha x)$ is an $\alpha$-contraction in the sup-norm, and denote by $\xi_\alpha \in \R^n$ the unique fixed point of this map. Then, we take as a certificate the germ at point $1^-$ of the function $\alpha \mapsto \xi_\alpha$. It follows from~\cite[Thm.~5.(b)]{Veinott1974} that this germ is represented by a Laurent series expansion of the form $a_{-1} / (1 - \alpha) + a_0 + a_1 (1 - \alpha) + \cdots$, with $a_{-1}, a_0, a_1, \ldots \in \mathbb{Q}^n$, and that this Laurent series is uniquely determined by the initial vectors $a_{-1}, \dots, a_{n - 1}$, which have short rational entries. The vector $a_{-1}$ is constant, all its coordinates being equal to $\rho$, recalling that the condition $\rho>0$ characterizes tropical independence. We have $T(\alpha \xi_\alpha) = \xi_\alpha$ for all $\alpha < 1$ close enough to $1$ if and only if the vectors $a_{-1}, \dots, a_{n - 1}$ satisfy a lexicographic system (Theorem~6, \emph{ibid.}), verifiable in polynomial time.
		\end{remark}
	
	\subsection{From stochastic mean-payoff games to tropical linear independence}
		
		We next show a converse to \Cref{thm:reduction}.
		
		\begin{theorem} \label{thm:reduction_converse}
			Checking whether the value of a stochastic turn-based mean-payoff game is nonpositive reduces in polynomial-time to checking whether a finite family of rational functions on a metric graph is tropically linearly dependent.
		\end{theorem}
		
		To establish this result, it will be convenient to consider the following constraint satisfaction problem. There are $n$ variables $c_1, \dots, c_n$. The input consists of a first collection of inequalities of the form
		\begin{align}
			c_i \geq \frac{c_j + c_k}{2}, \qquad (i, j, k) \in \mathcal{I}, \label{eq:CSP1}
		\end{align}
		for some subset $\mathcal{I} \subset [n]^3$ such that all triples $(i, j, k)$ have distinct elements and such that $j < k$; of a second collection of inequalities of the form
		\begin{align}
			c_i \geq \min(c_j, c_k), \qquad (i, j, k) \in \cJ, \label{eq:CSP2}
		\end{align}
		for some subset $\cJ \subset [n]^3$ such that again all triples $(i, j, k)$ have distinct elements and such that $j < k$; and of a third collection of inequalities
		\begin{align}
			c_i \geq a_{ij} + c_j, \qquad (i, j) \in \cK, \label{eq:CSP3}
		\end{align}
		where $\cK$ denotes the set of pairs $(i, j)$ of distinct elements of $[n]$, and where the scalars $a_{ij} \in \Z$ are given. We require that $a_{ij} + a_{ji} \leq 0$ holds for all such pairs; otherwise, the system is trivially unsatisfiable. Constraint satisfaction problems of this kind have been considered in~\cite{Bodirsky2017}.
		
		The next proposition follows from known results on turn-based games. We provide a proof for completeness.
		
		\begin{proposition} \label{prop:reduction_threshold}
			Checking whether the value of a stochastic turn-based mean-payoff game does not exceed a given threshold reduces in polynomial-time to checking whether a constraint satisfaction problem of the form~\eqref{eq:CSP1}--\eqref{eq:CSP3} admits a solution $c \in \R^n$.
		\end{proposition}
		
		\begin{proof}
			Andersson and Miltersen have shown in~\cite{andersson2009complexity} that both the quantitative problem (computing the value) and the strategic problem (computing optimal strategies) for general stochastic turn-based mean-payoff games reduce to the \emph{simple stochastic game} problem introduced by Condon~\cite{condon}, in which transition probabilities only take the values $0$, $1/2$ and $1$, and the payments are zero or one, being one only in an absorbing state.
			
			This builds on a reduction of Zwick and Paterson~\cite{zwick}, allowing one to replace constraints of the form $c_i \geq p \cdot c$, where $p$ is a rational probability vector, by systems of constraints involving probability vectors with half-integer or unit entries, as per~\eqref{eq:CSP1}, \eqref{eq:CSP2}, after the introduction of auxiliary variables.
			
			Moreover, Condon showed (Lemmas~7 and~8, \emph{ibid.}) that we always may assume that the \emph{value vector} $\bar{c} \in \R^S$ of a simple stochastic game with state space $S$ is the unique solution of a fixed-point equation of the type $c = F(c)$, where $F$ is an order-preserving self-map of $\R^S$, whose coordinates involve maxima, minima, and averages of variables as well as constant terms (p.~213, \emph{ibid.}). More precisely, the map $F$ is determined by a partition $S = S_{\max}\cup S_{\min} \cup S_{\textrm{avg}} \cup S_0 \cup S_1$, in which the sets $S_0$ and $S_1$ are non-empty, and by two maps $\alpha, \beta: (S_{\max}\cup S_{\min} \cup S_{\textrm{avg}}) \to S$, with
			\begin{align*}
				F_i(x) &= \diamond(x_{\alpha(i)}, x_{\beta(i)}), \qquad i \in S_{\diamond}, \; \diamond \in \{\max, \min, \operatorname{avg}\}, \\
				F_i(x) &= u, \qquad i \in S_u, \; u \in \{0,1\},
			\end{align*}
			in which $\operatorname{avg}$ denotes the function of two variables such that $\operatorname{avg}(u, v) = (u + v)/2$. We note that the requirement that the value $\bar{c}_{i^*}$ of a specified state $i^* \in S$ does not exceed a rational threshold $\theta$ is equivalent to the existence of a vector $c \in [0, 1]^S$ such that $c \geq F(c)$ and $\theta \geq c_{i^*}$. Indeed, if $\theta \geq \bar{c}_{i^*}$, the value vector $\bar{c}$ trivially satisfies these two conditions. Conversely, if a vector $c \in \R^S$ satisfies these conditions, using the order-preserving character of $F$, we deduce that $c \geq F(c) \geq \dots \geq F^k(c)$ for all $k \geq 1$. By uniqueness of the fixed-point of $F$, it follows that $\lim_{k} F^k(c) = \bar{c}$ and so $\theta_{i^*} \geq c_{i^*} \geq [F^k(c)]_{i^*} \geq \bar{c}_{i^*}$. Observe also that $\bar{c} \in [0, 1]^S$ since $F$ leaves the hypercube $[0, 1]^n$ invariant.
			
			We identify $S$ to $[n - 1]$, and introduce a dummy variable $c_n$, so that $c \in \R^S$ is now extended to a vector of $\R^n$. We next show that the system $[0, 1] \ni c_i \geq F_i(c), i \in S$, $\theta_{i^*} \geq c_{i^*}$, can be rewritten as~\eqref{eq:CSP1}-\eqref{eq:CSP3}. Indeed, we require that $c_n \leq c_i \leq c_n + 1$ for all $i \in S$ to encode the fact that $c_i \in[0, 1]$. We also rewrite the inequalities $c_{i^*} \leq \theta$ and $c_i \geq F_i(x) = 1$ for $i \in S_1$ as $c_{i^*} \leq \theta + c_n$ and $c_{i} \geq 1 + c_n$. Finally, we also add the redundant inequalities $c_i \geq c_j - 1$ for every pair $(i, j)$ of distinct elements of $[n]$. Multiplying $c$ by the denominator of the rational number $\theta$, we get that $[0, 1]^n \ni c \geq F(x)$ and $\theta_{i^*}\geq c_{i^*}$ is equivalent to a system of the form~\eqref{eq:CSP1}-\eqref{eq:CSP3}, in which every $a_{ij}$ is finite and integer valued.
		\end{proof}
		
		In order to show Theorem~\ref{thm:reduction_converse}, it will be convenient to consider the following \emph{generalized tropical linear dependence problem}. The input consists of metric graph $\Gamma$ \emph{not necessarily connected}, together with functions $f_1, \dots, f_n \colon \Gamma \to \R \cup \{\infty\}$, with $n \geq 2$, such that for every $i \in [n]$, the support of $f_i$, $\operatorname{supp}f_i \coloneqq \mleft\{x \in \Gamma \st f_i(x) < \infty\mright\}$ is a union of closed edges of the metric graph and of isolated vertices, and the restriction of $f_i$ to its support is continuous and piecewise-linear.
		
		We ask for the existence of real constants $c_1, \dots, c_n \in \R$ such that for all $x \in \Gamma$, the minimum $\min(c_1 + f_1(x), \dots, c_n + f_n(x))$ is achieved at least twice. Thus, by comparison with the tropical linear dependence problem considered in the first part of the paper, the generalization consists in relaxing the requirement that $\Gamma$ be connected and that the functions $f_i$ be everywhere finite.
		
		We shall derive~\Cref{thm:reduction_converse} from the following intermediate reduction.
		
		\begin{proposition} \label{lem:reduction_generalized_problem}
			The satisfiability problem for the system of constraints~\eqref{eq:CSP1}--\eqref{eq:CSP3} reduces in polynomial time to the generalized tropical linear dependence problem.
		\end{proposition}
		
		Our proof is based on the following construction. We associate to the collection of constraints~\eqref{eq:CSP1}--\eqref{eq:CSP3} the non-connected metric graph $\Gamma$ constructed as follows.
		
		We set $M \coloneq - \min_{ij} a_{ij} + 1$ and require that $M \geq 1$ (otherwise, $a_{ij} > 0$ would hold for all $(i, j) \in \cK$, contradicting our assumption that $a_{ij} + a_{ji} \leq 0$). Then,
		
		\begin{itemize}
			\item to every $(i, j, k) \in \cI$, we associate a closed edge $E_{ijk}$
			of the metric graph, of length $2M$, identified to the interval $[-M, M] \subset \R$;
			
			\item to every $(i, j, k) \in \cJ$, we associate two isolated vertices $v_{ijk}$ and $v'_{ijk}$;
			
			\item to every $(i, j) \in \cK$, we associate two other isolated vertices $w_{ij}$ and $w'_{ij}$.
		\end{itemize}
		
		Therefore, the metric graph $\Gamma$ is defined as the disjoint union of the edges $E_{ijk}$ as $(i, j, k) \in \cI$ and of the vertices $v_{ijk}$ and $v'_{ijk}$ with $(i, j, k) \in \cJ$ and $w_{ij}$ and $w'_{ij}$ with $(i, j) \in \cK$. In particular, each vertex and each of these edges constitutes a distinct connected component of $\Gamma$.
		
		We first define the tropical Dirac function at a point $y \in \Gamma$
		to be the function $\delta_y \colon \Gamma \to \R \cup \{\infty\}$, such that $\delta_y(x) = 0$ if $y = x$ and $\delta_y(x) = \infty$ otherwise. We now define a collection of functions $\Gamma \to \R \cup \{\infty\}$.
		
		\begin{itemize}
			\item For every $i \in [n]$, we define the function $f_i$ such that $f_i(x) = -x$ if $x$ belongs an an edge of the form $E_{lik}$, $f_i(x) = x$ if $x$ belongs to an edge of the form $E_{lji}$, and $f_i(x) = 0$ if $x$ belongs an an edge of the form $E_{ijk}$. Moreover, we set $f_i(v_{ljk}) = 0$ if $i \in \{l, j, k\}$,	$f_i(v'_{ljk}) = 0$ if $i \in \{j, k\}$, $f_i(w_{ij}) = 0$, $f_i(w_{ji}) = a_{ji}$, $f_i(w'_{ji}) = a_{ji}$, and we set $f_i(x) = \infty$ everywhere else.
			
			\item For every $(i, j, k) \in \cI$, we define two functions $f^\pm_{ijk}$ such that $f^\pm_{ijk}(x) = \mp x$ if $x \in E_{ijk}$, $f^\pm_{ijk}(x) = \infty$ everywhere else.
			
			\item For every $(i, j, k) \in \cJ$, we define the function $g_{ijk} = \min(\delta_{v_{ijk}}, \delta_{v'_{ijk}})$.
			
			\item For every $(i, j) \in \cK$, we also define the function $h_{ij} = \min(\delta_{w_{ij}}, \delta_{w'_{ij}})$.
		\end{itemize}
		
		The values taken by these functions are summarized in~\Cref{table:functions}, in which infinite values are omitted.
		
		{\renewcommand{\arraystretch}{1.4}
		\begin{table}
			\[ \begin{array}{c|ccccccccccccc}
				&
				E_{lik} & E_{lji} & E_{ijk} & v_{lik} & v_{lji} & v_{ijk} & v'_{lik} & v'_{lji} & v'_{ijk} & w_{ij} & w_{ji} & w'_{ij} & w'_{ji} \\
				\hline
				f_i
				& -x & x & 0 & 0 & 0 & 0 & 0 & 0 & & 0 & a_{ji} & & a_{ji} \\
				\hline
				f^\pm_{ijk}
				& & & \mp x \\
				\hline
				g_{ijk}
				& & & & & & 0 & & & 0 \\
				\hline
				h_{ij}
				& & & & & & & & & & 0 & & 0
			\end{array} \]
			\caption{The functions constructed in the proof of \Cref{lem:reduction_generalized_problem}.}
			\label{table:functions}
		\end{table}}
		
		\medskip
		
		We consider the following property.
		
		\begin{propertyD}
			For all $x \in \Gamma$, the minimum of the collection of terms
			\begin{equation}
				\begin{aligned}
					&c_i + f_i(x), i \in [n],
					&c^\pm_{ijk} + f^\pm_{ijk}(x), (i, j, k) \in \cI, \\
					&d_{ijk} + g_{ijk}(x), (i, j, k) \in \cJ,
					&d'_{ij} + h_{ij}(x), (i,j) \in \cK,
				\end{aligned}
				\label{eq:list}
			\end{equation}
			is achieved at least twice.
		\end{propertyD}

		\Cref{lem:reduction_generalized_problem} follows readily from the following lemma.
		
		\begin{lemma} \label{lem:equiv_csp_gendep}
			The vector $c \in \R^n$ satisfies the system of constraints~\eqref{eq:CSP1}--\eqref{eq:CSP3} if and only if Property~$\mathrm{D}$ holds for some choice of $c^\pm \in \R^{\cI}$, $d \in \R^{\cJ}$ and $d' \in \R^{\cK}$.
		\end{lemma}
		
		\begin{proof}
			\emph{``If'' part.} Suppose that Property~$\mathrm{D}$ holds. Taking $x = w_{ij}$, the minimum of the terms~\eqref{eq:list} specializes to
			\begin{align}
				\min \mleft(c_i, a_{ij} + c_j, d'_{ij}\mright).
				\label{eq:min-1}
			\end{align}
			
			To check this, the reader may find it convenient to refer to~\Cref{table:functions}. The different terms in the latter minimum arise from the values taken by the functions $f_i$, $f_j$ and $h_{ij}$ at point $w_{ij}$, which are visible in the columns $w_{ij}$ and $w_{ji}$.
			
			Similarly, taking $x = w'_{ij}$, the minimum of the terms~\eqref{eq:list} specializes to
			\begin{align}
				\min\mleft(a_{ij} + c_j, d'_{ij}\mright).
				\label{eq:min-2}
			\end{align}
			Since the minimum in~\eqref{eq:min-2} is achieved twice, we must have $a_{ij} + c_j = d'_{ij}$. Then, $c_i < a_{ij} + c_j$ would violate the assumption that the minimum in~\eqref{eq:min-1} is achieved at least twice. We deduce that $c_{i} \geq a_{ij} + c_j$ holds for all $(i, j) \in \cK$, which corresponds to the constraint given by~\eqref{eq:CSP3}.
			
			Taking $x \in E_{ijk}$, the minimum of the terms~\eqref{eq:list} specializes
			to
			\begin{align}
				\min\mleft(c_i, c_j - x, c_k + x, c_{ijk}^+ - x, c_{ijk}^- + x\mright)
				\label{eq:min-0}
			\end{align}
			If $c_i < \min(c_j - x, c_k + x)$, arguing by continuity, we get that $c_i < \min(c_j - y, c_k + y)$ holds at every point $y$ in a neighborhood of $x$. Since all the functions of $x$ arising in the minimum~\eqref{eq:min-0}, except the constant term $c_i$, have a derivative $\pm 1$, it is impossible that this minimum be achieved twice at every point of this neighborhood. It follows that $c_i \geq \min(c_j - x, c_k + x)$ holds for all $x \in E_{ijk}$.
			
			Applying Lemma~\ref{lem:dual_program} to compute the maximum of the latter minimum over the edge $E_{ijk}$, we deduce that $c_i \geq \min(c_j + M, c_k + M, (c_j + c_k)/2)$. We already proved that $c_i \geq a_{ij} + c_j$, which entails that $c_i > c_j - M$ for all $i \neq j$. It follows that $(c_j + c_k)/2 - (c_j + M) = (c_k - c_j)/2 - M < -M/2 < 0$. Hence, the minimum $\min(c_j + M, c_k + M, (c_j + c_k)/2)$ is achieved by the last term, and so, $c_i \geq (c_j + c_k)/2$ holds for all $(i, j, k) \in \mathcal{I}$, which is the constraint given by~\eqref{eq:CSP2}.
			
			Now, taking $x = v_{ijk}$, the minimum of the terms~\eqref{eq:list} specializes
			to
			\begin{align}
				\min(c_i, c_j, c_k, d_{ijk})
				\label{eq:min-3}
			\end{align}
			whereas taking $x = v'_{ijk}$, the same minimum specializes to
			\begin{align}
				\min(c_j, c_k, d_{ijk}).
				\label{eq:min-4}
			\end{align}
			If $c_i < \min(c_j, c_k)$, since the minimum in~\eqref{eq:min-3} is achieved twice, we must have $d_{ijk} = c_i$, but this contradicts the fact that the minimum in~\eqref{eq:min-4} is achieved twice. We deduce that $c_i \geq \min(c_j, c_k)$ holds for all $(i, j, k) \in \cJ$, which is~\eqref{eq:CSP1} and therefore, $c$ satisfies the collection of constraints~\eqref{eq:CSP1}--\eqref{eq:CSP3}.
			
			\medskip
			
			\emph{``Only if'' part.} Conversely, suppose that $c$ satisfies the collection of constraints~\eqref{eq:CSP1}--\eqref{eq:CSP3}. Since $c_i \geq a_{ij} + c_j$ holds for all $(i, j) \in \cK$, setting $d'_{ij} = a_{ij} + c_j$, we see that the minimum in each of the terms~\eqref{eq:min-1}, \eqref{eq:min-2} is attained at least twice.
			
			Similarly, since for all $(i, j, k) \in \cI$, $c_i \geq (c_j + c_k)/2$, setting $c_{ijk}^+ = c_j$ and $c_{ijk}^- = c_k$, we also deduce that the minimum in~\eqref{eq:min-0} is achieved twice at least at every point $x \in E_{ijk}$.
			
			Finally, since for all $(i, j, k) \in \cJ$, $c_i \geq \min(c_j, c_k)$, setting $d_{ijk} = \min(c_j, c_k)$, we see that the minimum in~\eqref{eq:min-3} and \eqref{eq:min-4} is achieved at least twice.
			
			Therefore, the vectors $c, c^\pm, d, d'$ satisfy Property~$\mathrm{D}$.
		\end{proof}
		
		To derive~\Cref{thm:reduction_converse} from \Cref{lem:reduction_generalized_problem}, we need the following lemma.
		
		\begin{lemma} \label{lem:bound}
 			Suppose that Property~$\mathrm{D}$ holds. Then, the vectors $c, c^\pm, d, d'$ involved in this property can always be chosen such that
			\begin{align} \label{eq:bound}
				c_i \in [0, M], \qquad c^\pm_{ijk} \in [0, M], \qquad d_{ijk} \in [0, M], \qquad d'_{ij} \in [-M, M].
			\end{align}
		\end{lemma}
		
		\begin{proof}
			Possibly after translating all entries of $c, c^\pm, d, d'$ by the same constant, we may assume that $0 = \min_{i} c_{i}$, in particular, $c_{i^*} = 0$ for some $i^*$. We saw in the proof of~\Cref{lem:reduction_generalized_problem} (paragraph after \eqref{eq:min-2}) that $c_i \geq a_{ij} + c_j$ for all pairs of distinct elements $(i,j )$, in particular $c_j \leq c_{i^*} - a_{i^*j} \leq M$ holds for all $j \in [n]$, which shows the first claim in~\eqref{eq:bound}.
			
			We also saw in this proof (paragraph after \eqref{eq:min-2}) that $d'_{ij} = a_{ij} + c_j \leq c_i$, which entails that $d'_{ij} \in [-M, M]$, showing the last claim in~\eqref{eq:bound}.
			
			Moreover, having fixed $c$, we note that we may always assume that $d_{ijk} \in \{c_i, c_j, c_k\}$. Otherwise, the term $d_{ijk}$ cannot attain the minima~\eqref{eq:min-3} and~\eqref{eq:min-4}, and we may always decrease $d_{ijk}$ until it meets one of these values (doing so trivially preserves the requirement that each of these minima be achieved at least twice). Therefore, we may assume that $d_{ijk} \in \{c_i, c_j, c_k\}$, and then $d_{ijk}$
			satisfies the third claim in~\eqref{eq:bound}.
			
			Finally, we observe that it is impossible that the minimum in~\eqref{eq:min-0} be achieved by $c_j - x$ for all $x \in [-M, M]$. Otherwise, we would have $c_j - x \leq c_k + x$ for all $x \in [-M, M]$, and taking $x = -M$, we get $2M < c_k - c_j \leq M$, a contradiction. Similarly, the same minimum cannot be achieved by $c_k + x$ for all $x \in [-M, M]$.
			
			We also saw in the proof of \Cref{lem:reduction_generalized_problem} (paragraph after \eqref{eq:min-0}) that we always have $c_i \geq \min(c_j - x, c_k + x)$ for all $x \in [-M, M]$. Let $E^\pm_{ijk} \coloneqq \{x \in E_{ijk} \mid c_j - x \lessgtr c_k + x\}$. We have $E_{ijk} = E^+_{ijk} \cup E^-_{ijk}$. We distinguish three cases.
			
			\emph{Case 1}. First assume that $E^+_{ijk} = E_{ijk}$. Then, the minimum in~\eqref{eq:min-0} coincides with $\min(c_j - x, c^+_{ijk} - x, c^-_{ijk} + x)$ for all $x \in E_{ijk}$, and since this minimum is achieved at least twice, we must have $c_j = c^+_{ijk}$. Moreover, replacing $c^-_{ijk}$ by $c_k$ preserves the property that this minimum is achieved at least twice at every point of $E_{ijk}$.
			
			\emph{Case 2}. Assume that $E^-_{ijk} = E_{ijk}$. We conclude by the dual argument that $c_k = c^-_{ijk}$ and that we may assume that $c_j = c^+_{ijk}$.
			
			\emph{Case 3}. Finally, assume that $E^+_{ijk}$ and $E^-_{ijk}$ are closed intervals with non-empty interior. Then, arguing as in Case~1, we see that we necessarily have $c_j = c^+_{ijk}$ and $c_k = c^-_{ijk}$.
			
			We showed that in all cases, we have, or may always assume that, $c_j = c^+_{ijk}$ and $c_k = c^-_{ijk}$. We deduce that the second claim in~\eqref{eq:bound} holds.
		\end{proof}
		
		We now complete the metric graph $\Gamma$ by adding one edge between every two vertices of $\Gamma$ that are not already connected, recalling that the vertices of $\Gamma$ are the endpoints of the edges $E_{ijk}$ together with isolated vertices. Each of these new edges has length two and will be identified with the real interval $[-1, 1]$. The completed metric graph $\barGamma$ obtained in this way is trivially connected.
		
		Moreover, we extend the functions constructed in the proof of~\Cref{lem:reduction_generalized_problem}. Let $\Phi$ denote the collection of all functions of the form $f_{i}, f_{ijk}^\pm, g_{ijk}$ or $h_{ij}$, and let $\varphi \in \Phi$. We first observe that:
		\begin{align}
			\text{for all }\varphi \in \Phi\text{ and } x \in \Gamma, \text{ either } \varphi(x) \in [-M, M] \text{ or } \varphi(x) = \infty.
			\label{eq:range_phi}
		\end{align}
		Indeed, when $\varphi = f^\pm_{ijk}= \mp x$, this follows from the fact that $|x| \leq M$ on every edge $E_{ikl}$; whereas when $\varphi = f_i$, the same properties follows by combining this fact with the inequalities $a_{ij} \geq -M + 1$ and $a_{ji} \leq -a_{ij} \leq M - 1$. Moreover,~\eqref{eq:range_phi} is trivial in the remaining cases.
		
		We next define the function $\barphi \colon \barGamma \to \R$, so that it coincides with $\varphi$ on $\supp \varphi$. The function $\barphi$ takes the value $4M$ on every edge or isolated vertex of $\Gamma$ on which $\varphi$ is infinite. So far, this defines a function $\barphi$ which is finite on every edge and vertex of the original metric graph $\Gamma$.
		
		Finally, we extend $\barphi$ to every added edge $E' \simeq [-1, 1]$, in such a way that
		\begin{align}
			\barphi(x) = \min(\barphi(-1) + 3M(x + 1), \barphi(1) - 3M(x - 1)). \label{eq:added}
		\end{align}
		Since $\barphi(-1), \barphi(1) \in [-M, 4M]$, the value of the latter minimum is achieved by the first term for $x = -1$. By symmetry, the value of the same minimum is achieved by the second term for $x = 1$. It follows that this expression of $\barphi$ provides a bona fide extension of $\varphi$ from $\{-1, 1\}$ to $[-1, 1]$. Moreover, the right derivative of $\barphi$ is $3M$ at point $-1$ whereas its left derivative is $-3M$ at point $1$.
		
		\medskip
		
		We next consider the following property.
		
		\begin{propertyDbar}
			For all $x$ in the completed metric graph $\barGamma$, the minimum of the collection of terms
			\begin{equation}
				\begin{aligned}
					&c_i + \overline f_i(x), i \in [n],
					&c^\pm_{ijk} + \overline f^\pm_{ijk}(x), (i, j, k) \in \cI, \\
					&d_{ijk} + \overline g_{ijk}(x), (i, j, k) \in \cJ,
					&d'_{ij} + \overline h_{ij}(x), (i, j) \in \cK,
				\end{aligned}
				\label{eq:list_bar}
			\end{equation}
			is achieved at least twice.
		\end{propertyDbar}
		
		We say that Property~$\mathrm{D}$ (or $\overline{\mathrm{D}}$) is \emph{satisfiable} if it is holds for some choice of vectors $c, c^\pm, d, d'$.
		
		\begin{lemma} \label{lem:gentoordinary}
			Property~$\mathrm{D}$ is satisfiable if and only if Property~$\overline{\mathrm{D}}$ is satisfiable.
		\end{lemma}
		
		\begin{proof}
			\emph{``Only if'' part.} Suppose that Property~$\mathrm{D}$ holds for some choice of vectors $c, c^\pm, d, d'$, which we require to lie in the intervals~\eqref{eq:bound}.
			
			Then, at every point $x$ of the original metric graph $\Gamma$, at least two terms are finite and achieve the minimum in~\eqref{eq:list}. Moreover, using the bounds~\eqref{eq:bound} and~\eqref{eq:range_phi}, we deduce that these two terms belong to $[-2M, 2M]$. Moreover, we defined $\barphi$ in such a way that $\barphi(x) = 4M$ for all $x \in \Gamma$ such that $\varphi(x) = \infty$. Since none of the scalars in~\eqref{eq:bound} is smaller than $-M$, we deduce that at such a point $x$, the term corresponding to $\barphi$ in~\eqref{eq:list_bar} has a value not lower than $3M$, and therefore it cannot contribute to the minimum of the terms~\eqref{eq:list_bar}. It follows that Property~$\overline{\mathrm{D}}$ holds for all $x\in \Gamma$.
			
			It remains to show that this property holds when $x$ belongs to an edge $E'$ added when constructing $\barGamma$ from $\Gamma$. On this edge, all the functions $\overline f_i$, $\overline f^\pm_{ijk}$, $\overline g_{ijk}$, $\overline h'_{ij}$ are defined as in~\eqref{eq:added}. In particular, identifying $E' \simeq [-1, 1]$, each of these functions has right derivative $3M$ at point $-1$ and left derivative $-3M$ at point $1$. We also know that the minimum in~\eqref{eq:list_bar} is achieved at least twice at every endpoint of the interval $E'$, since this endpoint belongs to the original graph $\Gamma$.
			
			Denoting by $C_{\pm 1}$ the value of the minimum of~\eqref{eq:list_bar} at $x = \pm 1$, we deduce that the minimum of~\eqref{eq:list_bar} coincides with $\psi(x) \coloneqq \min(C_{-1} + 3M(x + 1), C_1 - 3M(x - 1))$ for all $x \in E'$, and that the terms which achieve the minimum in ~\eqref{eq:list_bar} for $x = -1$ still achieve the minimum in~\eqref{eq:list_bar} at every point $x$ of $E'$ for which $\psi(x) = C_{-1} + 3M(x+1)$. A dual conclusion applies to $x = 1$. This implies that the minimum of the terms~\eqref{eq:list_bar} is achieved twice for all $x \in E'$, concluding the proof that Property~$\overline{\mathrm{D}}$ holds.
			
			\medskip
			
			\emph{``If'' part.}	Suppose that Property~$\overline{\mathrm{D}}$ holds for some choice of $e = (c, c^\pm, d, d')$. After translating all the entries of $e$ by a constant, we require the minimal entry of $e$ to be $0$. We denote by $\overline{m}(x)$ the minimum of the terms in~\eqref{eq:list_bar}.
			
			We first make the following observation, where $\alpha \geq 0$:
			\begin{align}
				(c_i \leq \alpha \leq 3M \text{ for some } i \in [n]) \implies (c_j < \alpha + M \text{ for all } j \in [n]). \label{e:transfer}
			\end{align}
			
			Indeed, if $c_i \leq \alpha$, evaluating the minimum of the terms~\eqref{eq:list_bar} at $x = w_{ij}$, we see that $4M > \alpha \geq c_i= c_i + f_i(w_{ij}) \geq \overline{m}(w_{ij})$, so that the terms in~\eqref{eq:list_bar} arising from the functions $\varphi \in \Phi$ such that $\barphi(w_{ij}) = 4M$ do not contribute to this minimum. In particular, $\overline{m}(w_{ij}) = \min(c_i, a_{ij} + c_j, d'_{ij})$ holds. Since the minimum $\overline{m}(w_{ij})$ is achieved at least twice, one of the following cases applies.
			
			\emph{Case~(i)}. $c_i = a_{ij} + c_j$. Then $c_j = -a_{ij} + c_i < \alpha + M$.
			
			\emph{Case~(ii)}. $d'_{ij} = c_i$. Then, we deduce by the same method that $\overline{m}(w'_{ij}) = \min(a_{ij} + c_j, d'_{ij})$, and since this minimum is achieved twice, $a_{ij} + c_j = d'_{ij} = c_i$ so that $c_j = -a_{ij} + c_i < \alpha + M$.
			
			\emph{Case~(iii)}. $d'_{ij} = a_{ij} + c_j < c_i$. Then $d'_{ij} < \alpha$. We argue again that $\overline{m}(w'_{ij}) = \min(a_{ij} + c_j, d'_{ij})$ is achieved twice and deduce as above that $c_j < \alpha + M$, showing~\eqref{e:transfer}.
			
			We now distinguish further cases, depending on the entry achieving the minimum in the vector $e$.
			
			\emph{Case (a)}. First assume that the minimum of the entries of $e$ is achieved by some $c_i$, so that $c_i = 0$. Then, we deduce from~\eqref{e:transfer} that $c_j < M$ holds for all $j \in [n]$, and considering~\eqref{eq:min-1}--\eqref{eq:min-4}, we deduce that at any point $x \in \Gamma$, a function $\varphi \in \Phi$ such that $\barphi(x) = 4M$ cannot contribute to the minimum $\overline{m}(x)$, so that Property~$\mathrm{D}$ holds. Furthermore, observe that the weaker assumption that $c_j < 3M$ for all $j \in [n]$ suffices to reach this conclusion.
			
			\emph{Case (b)}. Now assume that the same minimum is achieved by some $d'_{ij}$, so that $d'_{ij} = 0$. Then $\overline{m}(w'_{ij}) = \min(a_{ij} + c_j, d'_{ij})$, and since this minimum is achieved twice $c_j = -a_{ij} + d'_{ij} < M$. Then, by~\eqref{e:transfer}, we deduce that $c_i < 2M$, and we conclude as in the proof of Case (a) that Property~$\mathrm{D}$ holds.
			
			\emph{Case (c)}. Assume that the same minimum is achieved by some $d_{ijk}$, so that $d_{ijk} = 0$. Then, $\overline{m}(v_{ijk}) = \min(c_i, c_j, c_k, d_{ijk}) = 0$, and this minimum is achieved at least twice, which entails that at least one entry of $c$ is zero. So we are reduced to Case (a).
			
			\emph{Case (d)}. Finally, assume that the same minimum is achieved by some $c^\pm_{ijk}$, say by $c^+_{ijk}$, so that $c^+_{ijk} = 0$. Then, for all $x \in E_{ijk}$, we have $\overline{m}(x) = \min(c_i, c_j - x, c_k + x, c_{ijk}^+ - x, c_{ijk}^- + x)$, and since this minimum is achieved at least twice, we deduce that either $c_i \leq c^+_{ijk} + M = M $ or $c_j \leq c_{ijk}^+ = 0$, and therefore conclude as in Case (a).	
		\end{proof}
		
		We now conclude the proof of Theorem~\ref{thm:reduction_converse}. By \Cref{prop:reduction_threshold}, stochastic turn-based mean-payoff games reduce to the satisfiability of the constraint satisfaction problem~\eqref{eq:CSP1}--\eqref{eq:CSP3}, which, by~\Cref{lem:reduction_generalized_problem}, reduces to the generalized tropical linear dependence problem, which, in turn, by~\Cref{lem:gentoordinary}, reduces to the tropical linear dependence problem. \hfill \qed
	
	\subsection{Finite evaluation maps on semimodules of rational functions} \label{subsec:finite_evaluation}
		
		Given a subsemimodule $M$ of $R(D)$, it is natural to ask whether functions $f \in M$ can be fully characterized by the set of their values on a well-chosen finite set of points in $\mg$. More precisely:
		
		\begin{question} \label{q:finite_measurement_function}
			Given a subsemimodule $M$ of $R(D)$, do there exist a positive integer $K$ and points $x_1, \dots, x_K \in \mg$ such that the evaluation map $\eval \colon M \to \TT^K$, defined by
			\[ \eval(f) = (f(x_k))_{1 \leq k \leq K}, \]
			is injective?
		\end{question}
		
		This is motivated by the following observation.
		
		\begin{proposition} \label{prop:rank_projection}
			Suppose that $M$ is a module for which the answer to Question~\ref{q:finite_measurement_function} is positive. Then, $\troprank(M) = \troprank(\eval(M))$.
		\end{proposition}
		
		\begin{proof}
			A family $\eval(f_1), \dots, \eval(f_r)$ with $f_1, \dots, f_r \in M$ is tropically dependent if and only if there are real numbers $\lambda_1, \dots, \lambda_r$ such that for all $k \in [K]$, the minimum in $\min_{j \in [r]} (\lambda_j + f_j(x_k))$ is achieved at least twice. This condition can be rewritten as
			\[
				\eval(g) = \eval(g_j) \text{ for all $j$, with } g = \min_{s \in [r]} (\lambda_s + f_s) \text{ and } g_j = \min_{s \in [r] \setminus \{j\}} (\lambda_s + f_s).
			\]
			If $\eval$ is injective, if follows that $g = g_j$ holds for all $j \in [r]$. This means that $f_1, \dots, f_r$ are tropically dependent. We infer that $\troprank(M) \leq \troprank(\eval(M))$. The other inequality is trivial.
		\end{proof}
		
		The answer to Question~\ref{q:finite_measurement_function} depends on the structure of $M$, namely, on the number of distinct slopes that functions in $M$ realize along any edge of the metric graph. To formulate this dependence, we refine the combinatorial model of $\mg$ so that the set of slopes taken by functions $f \in M$ along unit tangent directions becomes constant on each edge. More precisely, for each pair $(e, v)$ consisting of an edge $e$ and an extremity $v$ of $e$, and any point $x$ of $\mg$ lying on the half-closed interval $e \setminus \{v\}$, let $\nu_x$ denote the unit tangent vector at $x$ directed toward $v$. Then, we require that the set $\Sl_{(e, v)}(M) \coloneqq \mleft\{\Sl_{\nu_x} f \st f \in M\mright\}$ be independent of $x$. A compactness argument, using the closedness of $R(D)$, ensures that such a model always exists.
		
		Let us first assume that $\Sl_{(e,v)}(M)$ has at most two elements for every pair $(e, v)$ in the graph. This is the case, for instance, if $\troprank(M) \leq 1$ (see Proposition~\ref{prop:comparison}). Since, on each edge $e$ with endpoints $u$ and $v$ and for every $f \in M$, the function $f$ changes slope at most once along $e$, the values of $f(u)$ and $f(v)$ completely determine $f$ on the entire edge $e$. Consequently, if $x_1, \dots, x_K$ are chosen as the vertices of $\mg$, the corresponding evaluation map $\eval$ is injective.
		
		On the contrary, let us assume that $\Sl_{(e, v)}(M)$ contains at least three elements for some pair $(e, v)$. In this case, for any choice of points $x_1, \dots, x_K \in \mg$, the associated evaluation map $\eval$ is generally not injective. This means that no finite evaluation map can fully characterize the functions $f \in M$.
		
		To illustrate this, we further refine the combinatorial model of $\mg$ and assume the following:
		\begin{enumerate}
			\item there exists an edge $e'$ (in this new model) contained in the original edge $e$ such that there are functions $f_1, f_2, f_3 \in M$ with constant slopes $s_1, s_2, s_3$ on $e'$ satisfying $s_1 < s_2 < s_3$; and
			
			\item no point $x_i$ lies in the interior of $e'$, i.e., $x_i \notin \mathring e'$ for all $i$.
		\end{enumerate}
		
\begin{figure}[h!]
			\centering
			\begin{minipage}{0.45\textwidth}
				\centering
				\begin{tikzpicture}[scale=.85]
					\draw[->] (0,0) -- (3.5,0);
					\draw[->] (0,0) -- (0,3.5);
					\draw[color=teal,semithick] (0,0) -- (3,0);
					\draw[color=red,semithick] (0,0) -- (3,1.5);
					\draw[color=blue, semithick] (0,0) -- (3,3);
					\draw[color=teal] (3,0) node[below right]{$f_1$};
					\draw[color=red] (3,1.5) node[right]{$f_2$};
					\draw[color=blue] (3,3) node[right]{$f_3$};
					\draw[color=teal] (0.8,0) node[below]{$s_1$};
					\draw[color=red] (1.5,0.75) node[below right]{$s_2$};
					\draw[color=blue] (0.8,0.7) node[above left]{$s_3$};
				\end{tikzpicture}
			\end{minipage}
			\begin{minipage}{0.45\textwidth}
				\centering
				\begin{tikzpicture}[scale=.85]
					\draw[->] (0,0) -- (3.5,0);
					\draw[->] (0,0) -- (0,3.5);
					\draw[color=blue,semithick] (0,0) -- (1.4,1.4);
					\draw[color=teal,semithick] (1.6,1.5) -- (3,1.5); % Before red lines for intersection colors
					\draw[color=red,semithick] (0.5,0.5) -- (2.5,1.5);
					\draw[color=red,semithick] (1,1) -- (2,1.5);
					\draw[color=red,semithick] (1.2,1.2) -- (1.8,1.5);
					\draw[color=red,semithick] (1.3,1.3) -- (1.7,1.5);
					\draw[color=red,semithick] (1.4,1.4) -- (1.6,1.5);
					\draw (3,1.5) node[right]{$f$};
					\draw[color=blue] (0.8,0.7) node[above left]{$s_3$};
					\draw[color=red] (1.5,1) node[below right]{$s_2$};
					\draw[color=teal] (2.25,1.5) node[above]{$s_1$};
					\phantom{\draw[color=teal] (3,0) node[below right]{$f_1$};} % For figure alignment
				\end{tikzpicture}
			\end{minipage}
			\caption{Construction of infinitely many functions using three different slopes.}
			\label{fig:three_slopes_distinct_functions}
		\end{figure}
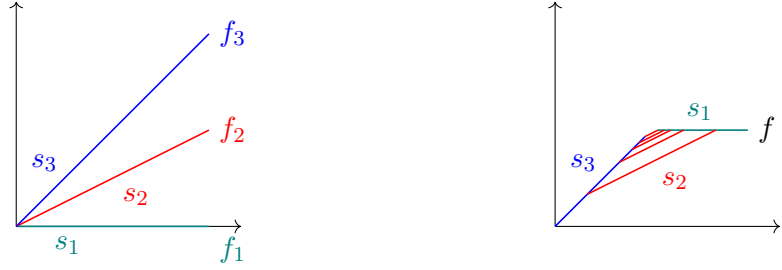

		It is easy to see that if $M$ contains sufficiently many functions -- for instance, if $M = R(D)$ -- then, using the functions $f_i$, we can construct infinitely many distinct functions $f \in M$ that coincide outside $\mathring e'$ (see Figure~\ref{fig:three_slopes_distinct_functions}). Since these functions share the same value at all the chosen points $x_i$, they have the same image under the evaluation map. Since this argument can be applied to any choice of the points $x_1, \dots, x_K$, for any $K$, this proves that the finite evaluation maps cannot be injective in general.
		
	\subsection{Geometric interpretation of the equivalence between turn-based games and tropical linear independence over metric graphs}
		
		Given a family of functions $f_1, \dots, f_r$ generating a submodule $M \subset \Rat(\mg)$, we define the \emph{ambiguity module} $\Amb(f_1, \dots, f_r) \subset M$ as follows. A function $h$ belongs to $\Amb(f_1, \dots, f_r)$ if either $h = \infty$, or there are real numbers $c_1, \dots, c_r$ such that $h(x) = \min_{1 \leq i \leq r} \, (f_i(x) + c_i) $ and the latter minimum is achieved at least twice for all $x \in \mg$.
		
		We say that a subset of $\R^K$ is \emph{semilinear} if it is a finite union of sets defined by finitely many weak or strict affine inequalities with rational coefficients.
		
		The following result provides a geometric interpretation of the equivalence established in~\Cref{thm:equiv}.
		
		\begin{theorem} \label{thm:expressive}
			Let $C \subset \R^K \cup \{\infty\}$ be any subset. The following assertions are equivalent:
			
			\begin{enumerate}[label=(A\arabic*)]
				\item \label{item:semilinear} $C$ is a submodule of $\mathbb{T}^K$, and the set $C \setminus \{\infty\}$ is closed, bounded with respect to Hilbert's seminorm, and semilinear.
				
				\item \label{item:ambiguity} There exist a metric graph $\mg$, functions $f_1, \dots, f_r \in \Rat(\mg) \setminus \{\infty\}$ and points $x_1, \dots, x_K$ of $\mg$ with rational coordinates such that $C = \vartheta(\Amb(f_1, \dots, f_r))$, where $\vartheta$ is the evaluation map on the points $x_1, \dots, x_K$.
			\end{enumerate}
		\end{theorem}
		
		Recall that the evaluation map $\vartheta$ on the points $x_1, \dots, x_K$ was defined in Question~\ref{q:finite_measurement_function}. Moreover, the term ``rational coordinates'' is understood recalling our assumption that the metric graph $\mg$ is given by a model with rational edge lengths; we require every point $x_i$ to be a barycenter with rational weights of the vertices of the edge to which it belongs.
		
		Sets $C$ satisfying the first condition of \Cref{thm:expressive} have been studied in~\cite{heltonnie}. In particular, it follows from Proposition~24, \emph{ibid.}, that $C \setminus \{\infty\}$ is a polyhedral complex whose cells have walls with normals of the form $\pm (e_i - p)$, where $p$ is a stochastic vector (nonnegative vector of sum one) and $e_i$ is the $i$-th basis vector. This contrasts with sets of the form $\vartheta(M)$ where $M$ is a finitely generated submodule of $\Rat(\mg)$, which turn out to be more special. Indeed, $\vartheta(M)$ is a ``tropical polyhedron'' in the sense of~\cite{DS04}; in particular, it is a polyhedral complex whose cells have walls with normals in the $A_n$ root system, i.e., normals of the form $e_i - e_j$ with $i \neq j$ (in other words, the stochastic vector $p = e_j$ is trivial).
		
		\begin{proof}[Proof of \Cref{thm:expressive}]
			\ref{item:semilinear} $\Rightarrow$ \ref{item:ambiguity}. By~\cite[Prop.~24]{heltonnie}, every set $C \subset \R^K \cup \{\infty\}$ satisfying \ref{item:semilinear} is such that there exists a Shapley operator $T$ of the form~\eqref{eq:def_Shapley}, with rational parameters $r_i^{\alpha \beta}$ and $P_{i,j}^{\alpha \beta}$ such that $C \setminus \{\infty\} = \{g \in \R^K \mid T(g) \leq g\}$. Moreover, Proposition~36, \emph{ibid}, implies that $C \setminus \{\infty\}$ can be rewritten as the projection of the set $S \subset \R^n$ of solutions of a constraint satisfaction problem of the form~\eqref{eq:CSP1}--\eqref{eq:CSP3}, for some $n \geq K$. This means that $C \setminus \{\infty\} = \pi(S)$, where $\pi$ is the map from $\R^n$ to $\R^K$ which projects on the first $K$ coordinates, $\pi(c_1, \dots, c_n) = (c_1, \dots, c_K)$. We now associate to the satisfaction problem~\eqref{eq:CSP1}--\eqref{eq:CSP3} a metric graph $\overline{\mg}$ and functions $\overline f_i$, $\overline f^\pm_{ijk}$, $\overline g_{ijk}$, $\overline h_{ij}$ as in the proof of \Cref{thm:reduction_converse} (paragraph preceding Property $\overline{\mathrm{D}}$), and denote by $M^{\operatorname{Amb}}$ the ambiguity module associated to this collection of functions. So, a map $h$ belongs to $M^{\operatorname{Amb}}$ if and only if for all $x \in \overline{\mg}$, $h(x)$ is the minimum of the collection of terms~\eqref{eq:list_bar}, and this minimum	is achieved at least twice. Then, we showed in the proof of~\Cref{thm:reduction_converse} that the minimum in~\eqref{eq:min-2} is still achieved twice, so that $h(w'_{ij}) = a_{ij} + c_j = d'_{ij}$. Then, taking the $K$ points $x_j \coloneqq w'_{1j}$ for $1 \leq j \leq K$ and considering the associated evaluation map $\vartheta$, we get that $\vartheta(M^{\operatorname{Amb}})$ is precisely the set of vectors $(a_{1j} + c_j)_{1 \leq j \leq K}$ where $c = (c_j)_{1 \leq j \leq K}$ belongs	to $\pi(S)$. Denoting by $\xi$ the vector $(a_{1j})_{1 \leq j \leq K}$, this shows that the translated semilinear set $\xi + C$ coincides with $\vartheta(M^{\Amb})$. Since the class of sets $C$ verifying~\ref{item:semilinear} is stable by translation, we deduce that any such set $C$ can be realized in the form~\ref{item:ambiguity}.
			
			\ref{item:ambiguity} $\Rightarrow$ \ref{item:semilinear}. Recall that semilinear sets coincide with sets that are definable in the first order theory of divisible ordered groups~\cite[Proof of Cor.~3.1.17]{Mar02}. It follows that $\vartheta(\Amb(f_1, \dots, f_r)) \setminus \{\infty\}$ is semilinear. Since the functions $f_1, \dots, f_r$ are bounded on $\mg$, all the vectors $\vartheta(f_i) \in \R^K$, with $1 \leq i \leq r$, are included in a ball in Hilbert's seminorm centered at zero and of a sufficiently large radius. Balls in Hilbert's seminorm are stable by tropical linear combinations, and thus, $\vartheta (\Amb(f_1, \dots, f_r))$ is bounded in Hilbert's seminorm. Moreover, since $\Amb(f_1, \dots, f_r)$ is a submodule of $\Rat(\mg)$, its image by the evaluation map is a submodule of $\mathbb{T}^K$. Finally, we observe that every $h \in \Amb(f_1, \dots, f_r)$ can be written as $h = \sup_{1 \leq i \leq r} (f_i + c_i)$, and that the real constants $c_1, \dots, c_K$ can always be chosen such that $c_i - c_j \leq \sup_{x \in \mg} (f_j(x) - f_i(x))$ for all $1 \leq i, j \leq K$ (otherwise, $c_i + f_j(x) > c_j + f_j(x)$ holds for all $x$ and the constant $c_j$ may be increased without changing $h$). By a compactness argument, this implies that $\vartheta(\Amb(f_1, \dots, f_K))$ is closed.
		\end{proof}

	\subsection{Computing the tropical rank}
			
		Computing the tropical rank of a matrix turns out to be a harder problem than checking tropical linear independence. In fact, this problem is known to be $\mathrm{NP}$-hard, even for matrices with entries in $\{0, 1\}$, see~\cite[Thm.~13]{kimroush}. We use this, together with the results of the previous section, to deduce the following analogue for rational functions.
		
		\begin{theorem} \label{thm:NP_hardness}
			Computing the tropical rank of finitely generated subsemimodules of rational functions on metric graphs is $\mathrm{NP}$-hard.
		\end{theorem}
		
		\begin{proof}
			Consider a matrix $A \in \TT^{m \times n}$ with entries in $\{0, 1\}$. First, we define a metric graph $\Gamma$ of model $(G = (V, E), \ell)$, with $G = (V, E)$ the complete graph on $m$ vertices and $\ell(e) = 2$ for all edges of $G$. Let $V = \{v_1, \dots, v_m\}$ be the vertex set and $E = \mleft\{\{v_i, v_s\} \st i, s \in [m], i \neq s\mright\}$ the edge set of $G$. For every edge $\{v_i, v_s\}$, denote by $w_{is}$ the midpoint of the edge.
			
			We associate to $A$ a semimodule $M$ of rational functions on $\Gamma$, generated by the following rational functions $f_1, \dots, f_n$. For each $j \in [n]$, set
			\[
				f_j(v_i) = A_{ij} \text{ and } f_j(w_{is})= \min(A_{ij}, A_{sj}) \text{ for all pair of distinct elements } i, s \in [m].
			\]
			We extend $f_j$ to $\Gamma$ by linear interpolation. It is easy to see that $f_j\in \Rat(\Gamma)$.
			
			Since the entries of the matrix $A$ belong to $\{0, 1\}$, for each edge $\{v_i, v_s\}$ of $G$, one of the following four cases occurs:
			\begin{align*}
				&\bullet \; f_j(v_i) = f_j(w_{ij}) = f_j(v_s) = 0, &&\bullet \; f_j(v_i) = f_j(w_{is}) = 0 < f_j(v_s) = 1, \\
				&\bullet \; f_j(w_{is}) = f_j(v_s) = 0 < f_j(v_i) = 1, &&\bullet \; f_j(v_i) = f_j(w_{ij}) = f_j(v_s) = 1.
			\end{align*}
			
			Now consider the family of points $x_1, \dots, x_K$ consisting of the vertices $v_1, \dots, v_m$ together will all the midpoints $w_{is}$ with $i < s$ and $i, s \in [m]$. We deduce by checking the above four cases that every $f_j$ satisfies the ``two-slopes'' condition stated in Subsection~\ref{subsec:finite_evaluation}, so that the restriction map $\eval \colon M \to \TT^{K}$ is injective.
			
			By Proposition~\ref{prop:rank_projection}, $\troprank(M) = \troprank(\eval(M)) = \troprank(B)$, where the matrix $B \in \TT^{K \times n}$ is given explicitly by $B_{kj} = A_{ij}$ if $x_k$ is equal to some vertex $v_i$, and $B_{kj} = \min(A_{ik}, A_{sj})$ if $x_k$ is equal to some midpoint $w_{is}$. Therefore, the rows of the matrix $B$ comprise all the rows of $A$ together with additional rows each of which is a tropical sum of pair of rows in $A$. It follows that $\troprank(B) = \troprank(A)$.
			
			By ~\cite[Thm.~13]{kimroush}, checking whether $\troprank(A) \geq r$ is an $\mathrm{NP}$-hard problem. We infer that checking whether $\troprank(M) \geq r$ for a finitely generated subsemimodule $M \subseteq \Rat(\Gamma)$ is also $\mathrm{NP}$-hard.
		\end{proof}

		\begin{remark}
			It follows from~\Cref{thm:equiv} that checking whether a subsemimodule $M$ of $\Rat(D)$ generated by $n$ functions $f_1, \dots, f_n$ has tropical rank at least $n - k$ reduces to solving ${n \choose k}$ stochastic turn-based mean-payoff games: it suffices to enumerate the ${n \choose k}$ subfamilies $f_{i_1}, \dots, f_{i_{n - k}}$ and to check which of these are tropically linearly independent. Similarly, checking whether $M$ has tropical rank at most $k$ reduces to solving ${n \choose {k + 1}}$ stochastic turn-based mean-payoff games. This is tractable for small values of $k$.
		\end{remark}

\section{Complementary results and further questions} \label{sec:complementary}
	
	In this section, we discuss several complementary results and questions, including a possible extension to higher dimension.
	
	\subsection{Closedness versus finite generation of subsemimodules}
		
		In general, $R(D)$ contains closed subsemimodules that are not finitely generated -- see the example below, which provides an answer to~\cite[Q.~6.2]{JP22}, although with a less rigid notion of tropical linear series. However, subsemimodules of $R(D)$ arising from the tropicalization of linear series on curves are always finitely generated, see~\cite[\S~9.4]{AG22}.
		
		\begin{example} \label{ex:infinite}
			Let $\Gamma$ be a metric graph with model $(G = (V, E), \ell)$, and let $x_1, x_2, x_3$ be three distinct points on an edge $e = \{u, v\}$ of $G$ in $\Gamma$, such that $x_2$ is the midpoint of the edge and lies at the midpoint of the segment joining $x_1$ and $x_3$. Consider the divisor $D = n \, (u) + n \, (v)$ for a sufficiently large positive integer $n$, and define $M \subset R(D)$ to be the set of all functions $f \in R(D)$ satisfying the inequality $-\varepsilon + 2 f(x_2) \geq f(x_1) + f(x_3)$, for $\varepsilon > 0$ small. Then, $M$ forms a closed subsemimodule of $R(D)$. However, for $\varepsilon > 0$ small enough, $M$ is not finitely generated.
			
			To see this, consider the evaluation map
			\[
				\eval \colon R(D) \to \TT^3, \qquad \eval(f) = (f(x_1), f(x_2), f(x_3)).
			\]
			If $M$ were finitely generated, then $\eval(M)$ would also be finitely generated, implying that it would have only finitely many extreme points.
			
			\begin{figure}[htbp]
				\centering
				\begin{tikzpicture}
					% Draw x-axis
					\draw (-3,0) -- (3,0);
					
					\draw[dotted] (-3,-1.75) -- (3,-1.75);
					\draw[dotted] (-3,-0.75) -- (3,-0.75);
					\draw[dotted] (-3,-0.25) -- (3,-0.25);
					\node[left] at (-3, -0.25) {\footnotesize $-\alpha$};
					\node[left] at (-3, -0.75) {\footnotesize $-\beta$};
					\draw[dotted] (-1,0) -- (-1,-1.75);
					\draw[dotted] (1,0) -- (1,-1.75);
					
					% Define points
					\node[above] at (-2,0) {\small $u$};
					\node[above] at (-1,0) {\small $x_1$};
					\node[above] at (0,0) {\small $x_2$};
					\node[above] at (1,0) {\small $x_3$};
					\node[above] at (2,0) {\small $v$};
					
					\node[below] at (1.6,-1.2) {\small $-n$};
					\node[above] at (1.4,-0.65) {\small $-1$};
					\node[below] at (-1.6,-0.9) {\small $1$};
					
					% Draw vertical markers
					\foreach \x in {-2,-1,0,1,2}
					\draw (\x,0.1) -- (\x,-0.1);
					
					% Polygonal curve
					\draw[blue,thick] (-0.25, 0) -- (0.75, 0);
					\draw[blue,thick] (0.75, 0) -- (0.75+1,-1);
					\draw[blue, thick] (-2, -1.75) -- (-0.25,0);
					\draw[blue, thick] (2, -1.75) -- (1.75, -1); 
					\draw[blue, thick] (-2, -1.75) -- (-2.5, -1.75);
					\draw[blue, thick] (2, -1.75) -- (2.5, -1.75);
					
					\draw (5.25,0) -- (7,0);
					\draw (7,0) -- (7,-1.75);
					\filldraw[fill=lightgray, draw=black] (6,0) -- (6.75,0) -- (7,-0.25) -- (7,-1) -- (6,-1) -- cycle;
					\draw[red,very thick] (6.75,0) -- (7,-0.25);
				\end{tikzpicture}
				\caption{A parametric family of functions $f \in R(D)$ (left). Fragment of the cross section $\mleft\{(f(x_1), f(x_3)) \st f \in M, f(x_2) = 0\mright\}$ with an infinite upper Pareto set depicted in red (right).}
				\label{fig:closed_not_fg}
			\end{figure}
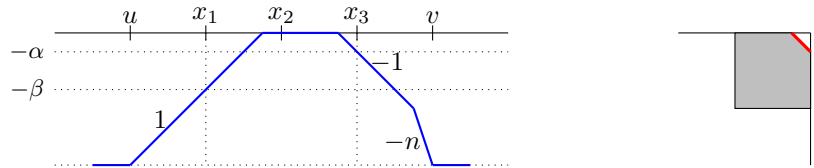
			
			Now, we shall use the following characterization of tropical extreme points, which can be found in~\cite{GK07} (proof of Theorem~3.1): if $N$ is a closed subsemimodule of $\TT^n$, then a point $z \in N$ is extreme if and only if there is an index $i \in [n]$ such that $z - z_i \, e_i$ is a maximal element of the cross section $S_i(N)$ of $N$ given by $S_i(N) = \mleft\{w \in N \st w_i = 0\mright\}$. Coming back to our example, it suffices to check that $S_2(\eval(M))$ has an infinite set of maximal elements.
			
			To show this, consider the function $f$ in $\Rat(\Gamma)$, constant on $\Gamma \setminus e$, and given on the edge $e$ by
			\[
				f(x) = \max(u - x_1 - \beta, \min(x - x_1 - \beta, x_3 - \alpha - x, u - x_1 - \beta + n (v - x))) \quad \text{for all $x$ in $e$},
			\]
			with $x$ a parameter on $e$ and $\alpha, \beta$ real numbers satisfying $0 \leq \alpha \leq \beta \leq L \coloneqq x_3 - x_2$. For $n$ large enough, this function has the shape shown in~Fig.~\ref{fig:closed_not_fg} (the slopes are indicated on the figure), with $\ord_u(f) = -1$, $\ord_v(f) = -n$, and $\ord_x(f) \geq 0$ elsewhere, thus $f \in R(D)$. Observe that $f(x_1) = - \beta$ and $f(x_3) = -\alpha$, where all the values $0 \leq \alpha \leq \beta \leq L$ are realizable.
			
			This entails that the set $\widetilde S_2 = \mleft\{(y_1, y_3) \st (y_1, 0, y_3) \in \eval(R(D))\mright\}$ contains every element of the form $(-\beta, -\alpha)$. By symmetry, we deduce that $\widetilde{S}_2 \supset [-L, 0]^2$. It follows that
			\[
				S_2(\eval(M)) \cap[-L,0]^2 = \mleft\{(y_1, y_3) \st (y_1 + y_3)/2 \leq -\varepsilon \text{ and } -L \leq y_1, y_2 \leq 0\mright\}.
			\]
			Hence, $S_2(\eval(M)) \cap [-L,0]^2$ has an infinite set of maximal elements, consisting of its upper boundary, defined by the equality $(y_1 + y_3) / 2 = -\varepsilon$. All the elements of this boundary are maximal elements of $S_2(\eval(M))$ as well, showing that the set of extreme points of $\eval(M)$ is infinite. This is illustrated in Figure~\ref{fig:closed_not_fg} (right).
		\end{example}
		
		When $M$ is finitely generated, by Theorem~\ref{thm:polyhedral_structure}, the linear system $|(D, M)|$ admits a polyhedral structure. In Section~\ref{sec:introduction}, for a more general subsemimodule $M$, we defined the dimension of the corresponding linear system as a supremum over finitely generated subsemimodules because as we show below (see Proposition~\ref{prop:o-minimal}), in general, $M$ can have pathological behavior, and we do not know the answer to the following question.
		
		\begin{question}
			Let $M \subset R(D)$ be closed. Is it possible to describe the geometric structure of $|(D, M)|$? In particular, does it admit a good notion of dimension?
		\end{question}
		
		Our next result shows the subtlety behind this question. Recall that o-minimal structures capture a notion of ``tameness'' for geometric objects; e.g., semialgebraic sets are the simplest examples of sets definable in an o-minimal structure. See~\cite{Dries98} for background.
		
		\begin{proposition} \label{prop:o-minimal}
			There exists a closed subsemimodule $M \subset R(D)$ such that $|(D,M)|$ is not definable in any o-minimal structure.
		\end{proposition}
		
		This result sheds light on the definability questions asked by Jensen and Payne~\cite[Q.~6.3]{JP22}.

		\begin{proof}
			We consider the metric graph $\Gamma$ with distinguished points $x_1, x_2, x_3$, and the divisor $D$ defined in Example~\ref{ex:infinite}. Let $\eval \colon R(D) \to \TT^3$ be the evaluation map on these three points.
			
			Let $\kappa \colon \R \to \R$ be an absolutely continuous function such that $0 \leq \kappa'(x) \leq 1$ holds almost everywhere. Now define
			the function $\psi \colon \R^2 \to \R$ by setting $\psi(y, z) = y + \kappa(z - y)$. Let $M_\kappa$ be the union of $\{\infty\}$ and the set of finite-valued functions $f \in R(D)$ which satisfy the inequality $f(x_2) \geq \psi(f(x_1), f(x_3))$. Given that $ 0 \leq \kappa'(x) \leq 1 $, it follows that
			\[
				\kappa(y) \leq \kappa(y + z) \leq \kappa(y) + z \quad \text{for all } y \in \mathbb{R}, \; z \in \mathbb{R}_{\geq 0}.
			\]
			This implies that $\psi$ is order-preserving. Moreover, it commutes with the addition of constants. Consequently, $ M_\kappa $ is a closed subsemimodule of $ R(D) $. For instance, the semimodule $M$ of Example~\ref{ex:infinite} corresponds to the special choice $\kappa(y) = \varepsilon / 2 + y / 2$.
			
			Now, we choose $\kappa$ so that $\kappa'(x)$ takes only the values $0$ and $1$ almost everywhere, and makes an infinite number of switches between the values $0$ and $1$ in a neighborhood of $0$. To make this example concrete, we may choose a sequence of positive real numbers $(\alpha_k)_k$ decreasing to zero, define the open interval $I_{k} = (\alpha_{k + 1}, \alpha_k)$, and set $\kappa'(x) = 0$ for all $x \in I_{2k} \cup (-I_{2k})$ and $\kappa'(x) = 1$ for all $x \in I_{2k + 1} \cup (-I_{2k + 1})$. Here, $-I_k$ is the interval $(-\alpha_k, -\alpha_{k + 1})$. We fix $\kappa(0) = -\varepsilon$ with $\varepsilon > 0$ small enough. Then, the cross section
			\[ S_2(\eval(M_\kappa)) = \{(f(x_1), f(x_2)) \mid f \in M_\kappa, f(x_2) = 0\} \]
			is a set similar to the one depicted in Figure~\ref{fig:closed_not_fg} (right), but now with an upper right boundary consisting of a staircase with an infinite number of stairs accumulating at the point $(-\varepsilon, -\varepsilon)$ of $\R^2$ with coordinates $-\varepsilon$. This implies that $S_2(\eval(M_\kappa))$ is not definable in any o-minimal structure.
			
			We claim that $|(D, M_\kappa)|$ is not a definable topological space in any o-minimal structure. Suppose for the sake of a contradiction that $|(D, M_\kappa)|$ is definable in some o-minimal structure.
			
			First, we embed $|(D, M_\kappa)|$ in $M_\kappa$ by sending each $E \in |(D, M_\kappa)|$ to the unique function $f \in M_\kappa$ with $E = D + \div(f)$ and $f(x_2) = 0$. The evaluation map $\eval \colon M_\kappa \to \TT^3$ defines a map $|(D, M_\kappa)| \to S_2(\eval(M_\kappa))$, and identifies $S_2(\eval(M_\kappa))$ with a topological quotient $|(D, M_\kappa)| / \mathcal{R}$ for a definable relation $\mathcal R \subset |(D, M_\kappa)| \times |(D, M_\kappa)|$. Note that $\mathcal R$ is definably proper in the sense that the preimage of any compact subset of $S_2(\eval(M_\kappa))$ under the projection map $|(D, M_\kappa)| \to S_2(\eval(M_\kappa))$ is compact. A theorem of van den Dries~\cite[Chap.~10, Thm.~2.15]{Dries98} then implies that $S_2(\eval(M_\kappa))$ is itself definable in the same o-minimal structure, leading to a contradiction.
		\end{proof}
	
	\subsection{Higher dimension}
		
		Let $Y \subseteq \R^d$ be a polyhedral subspace (e.g., a tropical subvariety). Let $\Rat(Y, \R)$ be the union of $\infty$ and the set of piecewise linear functions on $Y$ (with non-necessarily integral slopes). Endowed with the operation of tropical addition and tropical multiplication by constants, $\Rat(Y, \R)$ is a semimodule over $\TT$. Let $M$ be a finitely generated subsemimodule of $\Rat(Y, \R)$. For example, if $Y$ is a tropical subvariety, and $D$ is a divisor on $Y$, then $M$ may be a finitely generated subsemimodule of $R(D)$, where $R(D)$ is the union of $\{\infty\}$ and the set of piecewise linear functions on $Y$ with integral slopes. We define $\troprank(M)$ as the maximum integer $r$ such that there exist tropically independent elements $f_1, \dots, f_r \in M$.
		
		Let $g_1, \dots, g_l$ be a generating set for $M$. Consider the map 
		\[
			\Psi \colon \R^l \to M, \qquad (c_1, \dots, c_l) \mapsto \min_{j \in [l]} \bigl(g_j + c_j\bigr).
		\]
		We define a notion of dimension for $M$ as follows. Consider an element $f\in M$.
		
		\begin{proposition}
			The subset $\Psi^{-1}(f) \subset \R^l$ is polyhedral.
		\end{proposition}
		
		\begin{proof}
			Choose a polyhedral structure $\Delta$ on $Y$ such that the generators $g_1, \dots, g_l$ and $f$ are affine on each face of $\Delta$. For each face $\sigma \in \Delta$ and $(c_1, \dots, c_l) \in \Psi^{-1}(f)$, since $\min_{j \in [l]}(g_j + c_j)_{|\sigma} = f_{|\sigma}$, and $f_{|\sigma}$ is affine, we get the existence of $j = \mu(\sigma) \in [l]$ such that
			\[
				f_{|\sigma} = (g_{\mu(\sigma)} + c_{\mu(\sigma)})_{|\sigma} \leq (g_i + c_i)_{|\sigma}, \text{ for all } i \in [l].
			\]
			For each function $\mu \colon \Delta \to [l]$, let $C_\mu$ be the set of all points $(c_1, \dots, c_l) \in \R^l$ such that the inequality above is satisfied for all $\sigma \in \Delta$. This is a polyhedral subset of $\R^l$. Moreover, $\Psi^{-1}(f)$ is the union of the sets $C_\mu$. We infer the result.
		\end{proof}
		
		We define the dimension of $M$, denoted by $\dim(M)$, as
		\[
			\dim(M) = \max_{f \in M} \bigl[l - \dim \bigl(\Psi^{-1}(f)\bigr)\bigr].
		\]
		
		\begin{question}
			Let $M$ be a finitely generated subsemimodule of $\Rat(Y, \R)$. Do we have the equality $\troprank(M) = \dim(M)$?
		\end{question}
		
		A positive answer would generalize Theorem~\ref{thm:tropical_rank_equals_topological_rank}. Note that we may reduce to the case where $Y$ is compact in order to use Theorem~\ref{thm:certificate_independence_generalized}.
		
		We note that the stochastic game interpretation of the tropical linear independence problem carries over to the multidimensional case. However, by extending \Cref{rk:dual}, we see that the number of actions of the players is no longer polynomially bounded in the input size.

\bibliographystyle{plain}
\bibliography{Bibliography}

\begin{thebibliography}{10}

\bibitem{AGGut10}
Marianne Akian, St\'{e}phane Gaubert, and Alexander Guterman.
\newblock Tropical polyhedra are equivalent to mean payoff games.
\newblock {\em Internat. J. Algebra Comput.}, 22(1):1250001, 43, 2012.

\bibitem{akian2023}
Marianne Akian, St\'{e}phane Gaubert, Yang Qi, and Omar Saadi.
\newblock Tropical linear regression and mean payoff games: or, how to measure
  the distance to equilibria.
\newblock {\em SIAM J. Discrete Math.}, 37(2):632--674, 2023.

\bibitem{heltonnie}
Xavier Allamigeon, St\'{e}phane Gaubert, and Mateusz Skomra.
\newblock The tropical analogue of the {H}elton-{N}ie conjecture is true.
\newblock {\em J. Symbolic Comput.}, 91:129--148, 2019.

\bibitem{amini2013reduced}
Omid Amini.
\newblock Reduced divisors and embeddings of tropical curves.
\newblock {\em Trans. Amer. Math. Soc.}, 365(9):4851--4880, 2013.

\bibitem{AG22}
Omid Amini and Lucas Gierczak.
\newblock Limit linear series: combinatorial theory.
\newblock {\em J. Comb. Algebra}, 2026.
\newblock published online first.

\bibitem{andersson2009complexity}
Daniel Andersson and Peter~Bro Miltersen.
\newblock The complexity of solving stochastic games on graphs.
\newblock In {\em Algorithms and computation}, volume 5878 of {\em Lecture
  Notes in Comput. Sci.}, pages 112--121. Springer, Berlin, 2009.

\bibitem{BJ16}
Matthew Baker and David Jensen.
\newblock Degeneration of linear series from the tropical point of view and
  applications.
\newblock In {\em Nonarchimedean and tropical geometry}, Simons Symp., pages
  365--433. Springer, [Cham], 2016.

\bibitem{baker2007riemann}
Matthew Baker and Serguei Norine.
\newblock Riemann--{R}och and {A}bel--{J}acobi theory on a finite graph.
\newblock {\em Adv. Math.}, 215(2):766--788, 2007.

\bibitem{bewley_kohlberg}
Truman Bewley and Elon Kohlberg.
\newblock The asymptotic theory of stochastic games.
\newblock {\em Math. Oper. Res.}, 1(3):197--208, 1976.

\bibitem{Bodirsky2017}
Manuel Bodirsky and Marcello Mamino.
\newblock Tropically convex constraint satisfaction.
\newblock {\em Theory of Computing Systems}, 62(3):481--509, April 2017.

\bibitem{Butkovi2010}
Peter Butkovi\v{c}.
\newblock {\em Max-linear systems: theory and algorithms}.
\newblock Springer Monographs in Mathematics. Springer-Verlag London, Ltd.,
  London, 2010.

\bibitem{Butkovi1985}
Peter Butkovi\v{c} and Ferdinand Hevery.
\newblock A condition for the strong regularity of matrices in the minimax
  algebra.
\newblock {\em Discrete Appl. Math.}, 11(3):209--222, 1985.

\bibitem{Changetal}
Chih-Wei Chang, Matthew Dupraz, Hernan Iriarte, David Jensen, Dagan Karp, Sam
  Payne, and Jidong Wang.
\newblock Tropical linear series and matroids.
\newblock {\em Preprint arXiv:2508.20062}, 2025.

\bibitem{Chatterjee2014}
Krishnendu Chatterjee and Rasmus Ibsen-Jensen.
\newblock The complexity of ergodic mean-payoff games.
\newblock In {\em Automata, languages, and programming. {P}art {II}}, volume
  8573 of {\em Lecture Notes in Comput. Sci.}, pages 122--133. Springer,
  Heidelberg, 2014.

\bibitem{condon}
Anne Condon.
\newblock The complexity of stochastic games.
\newblock {\em Inform. and Comput.}, 96(2):203--224, 1992.

\bibitem{crandalltartar}
Michael~G. Crandall and Luc Tartar.
\newblock Some relations between nonexpansive and order preserving mappings.
\newblock {\em Proc. Amer. Math. Soc.}, 78(3):385--390, 1980.

\bibitem{develin2007rank}
Mike Develin, Francisco Santos, and Bernd Sturmfels.
\newblock On the rank of a tropical matrix.
\newblock In {\em Combinatorial and computational geometry}, volume~52 of {\em
  Math. Sci. Res. Inst. Publ.}, pages 213--242. Cambridge Univ. Press,
  Cambridge, 2005.

\bibitem{DS04}
Mike Develin and Bernd Sturmfels.
\newblock Tropical convexity.
\newblock {\em Doc. Math.}, 9:1--27, 2004.

\bibitem{Dupraz24}
Matthew Dupraz.
\newblock Tropical linear systems and the realizability problem.
\newblock {\em Master's thesis, École polytechnique fédérale de Lausanne,
  2024. arXiv:2506.21268}, 2024.

\bibitem{FJP20}
Gavril Farkas, David Jensen, and Sam Payne.
\newblock The {K}odaira dimensions of $\mathcal{M}_{22}$ and
  $\mathcal{M}_{23}$.
\newblock {\em Preprint arXiv:2005.00622}, 2020.

\bibitem{PFT}
St\'{e}phane Gaubert and Jeremy Gunawardena.
\newblock The {P}erron-{F}robenius theorem for homogeneous, monotone functions.
\newblock {\em Trans. Amer. Math. Soc.}, 356(12):4931--4950, 2004.

\bibitem{GK07}
St\'{e}phane Gaubert and Ricardo~D. Katz.
\newblock The {M}inkowski theorem for max-plus convex sets.
\newblock {\em Linear Algebra Appl.}, 421(2-3):356--369, 2007.

\bibitem{podolskii}
Dima Grigoriev and Vladimir~V. Podolskii.
\newblock Complexity of tropical and min-plus linear prevarieties.
\newblock {\em Comput. Complexity}, 24(1):31--64, 2015.

\bibitem{Groetschel1993}
Martin Gr\"{o}tschel, L\'{a}szl\'{o} Lov\'{a}sz, and Alexander Schrijver.
\newblock {\em Geometric algorithms and combinatorial optimization}, volume~2
  of {\em Algorithms and Combinatorics}.
\newblock Springer-Verlag, Berlin, second edition, 1993.

\bibitem{gurvich}
Vladimir~A. Gurvich, Alexander~V. Karzanov, and Leonid~G. Khachiyan.
\newblock Cyclic games and finding minimax mean cycles in digraphs.
\newblock {\em Zh. Vychisl. Mat. i Mat. Fiz.}, 28(9):1407--1417, 1439, 1988.

\bibitem{haase2012linear}
Christian Haase, Gregg Musiker, and Josephine Yu.
\newblock Linear systems on tropical curves.
\newblock {\em Math. Z.}, 270(3-4):1111--1140, 2012.

\bibitem{IRowen}
Zur Izhakian and Louis Rowen.
\newblock The tropical rank of a tropical matrix.
\newblock {\em Comm. Algebra}, 37(11):3912--3927, 2009.

\bibitem{JP14}
David Jensen and Sam Payne.
\newblock Tropical independence {I}: {S}hapes of divisors and a proof of the
  {G}ieseker-{P}etri theorem.
\newblock {\em Algebra \& Number Theory}, 8(9):2043--2066, 2014.

\bibitem{JP21}
David Jensen and Sam Payne.
\newblock Recent developments in {B}rill--{N}oether theory.
\newblock {\em Preprint arXiv:2111.00351}, 2021.

\bibitem{JP22}
David Jensen and Sam Payne.
\newblock Tropical linear series and tropical independence.
\newblock {\em Preprint arXiv:2209.15478}, 2022.

\bibitem{kimroush}
Ki~H. Kim and Fred~W. Roush.
\newblock Factorization of polynomials in one variable over the tropical
  semiring.
\newblock {\em Preprint arXiv:math/0501167v2}, 2005.

\bibitem{Mar02}
David Marker.
\newblock {\em Model theory}, volume 217 of {\em Graduate Texts in
  Mathematics}.
\newblock Springer-Verlag, New York, 2002.
\newblock An introduction.

\bibitem{sorin_repeated_games}
Jean-Fran\c{c}ois Mertens, Sylvain Sorin, and Shmuel Zamir.
\newblock {\em Repeated games}, volume~55 of {\em Econometric Society
  Monographs}.
\newblock Cambridge University Press, New York, 2015.
\newblock With a foreword by Robert J. Aumann.

\bibitem{shapley_stochastic}
Lloyd~S. Shapley.
\newblock Stochastic games.
\newblock {\em Proc. Nat. Acad. Sci. U.S.A.}, 39:1095--1100, 1953.

\bibitem{Dries98}
Lou van~den Dries.
\newblock {\em Tame topology and o-minimal structures}, volume 248 of {\em
  London Mathematical Society Lecture Note Series}.
\newblock Cambridge University Press, Cambridge, 1998.

\bibitem{Veinott1974}
Arthur~F. Veinott, Jr.
\newblock Markov decision chains.
\newblock In {\em Studies in optimization}, MAA Stud. Math., Vol. 10, pages
  124--159. Math. Assoc. America, Washington, DC, 1974.

\bibitem{yoshitomi}
Shuhei Yoshitomi.
\newblock Generators of modules in tropical geometry.
\newblock {\em Preprint arXiv:1001.0448}, 2010.

\bibitem{zwick}
Uri Zwick and Mike Paterson.
\newblock The complexity of mean payoff games on graphs.
\newblock {\em Theoret. Comput. Sci.}, 158(1-2):343--359, 1996.

\end{thebibliography}

\end{document}